\documentclass[leqno,11pt]{amsart}

\usepackage{amssymb, amsmath}
\def\norm#1#2{\|#1\|_{#2}}

\def\refer#1{~\ref{#1}}
\def\refeq#1{~(\ref{#1})}
\def\ccite#1{~\cite{#1}}

\def\longformule#1#2{
\displaylines{ \qquad{#1} \hfill\cr \hfill {#2} \qquad\cr } }
\def\inte#1{
\displaystyle\mathop{#1\kern0pt}^\circ }

\def\infetage#1#2{
\inf_{\scriptstyle {#1}\atop\scriptstyle {#2}} }

\def\supetage#1#2{
\sup_{\scriptstyle {#1}\atop\scriptstyle {#2}} }

\let\al=\alpha

\let\e=\varepsilon

\let\lam=\lambda

\let\s=\sigma
\let\f=\phi
\let\vf=\varphi

\let\D=\Delta
\let\Lam=\Lambda

\let\wt=\widetilde
\let\wh=\widehat

\def\cB{{\mathcal B}}
\def\cC{{\mathcal C}}
\def\cD{{\mathcal D}}
\def\cE{{\mathcal E}}
\def\cF{{\mathcal F}}

\def\cK{{\mathcal K}}
\def\cL{{\mathcal L}}

\def\cQ{{\mathcal Q}}

\def\cS{{\mathcal S}}

\def\tlk{T_{\Lambda,X}}

\def\virgp{\raise 2pt\hbox{,}}
\def\cdotpv{\raise 2pt\hbox{;}}

\def\eqdefa{\buildrel\hbox{{\rm \footnotesize def}}\over =}

\def\C{\mathop{\bf C\kern 0pt}\nolimits}
\def\DD{\mathop{\bf D\kern 0pt}\nolimits}
\def\K{\mathop{\bf K\kern 0pt}\nolimits}
\def\N{\mathop{\bf N\kern 0pt}\nolimits}
\def\Q{\mathop{\bf Q\kern 0pt}\nolimits}
\def\R{\mathop{\bf R\kern 0pt}\nolimits}
\def\SS{\mathop{\bf S\kern 0pt}\nolimits}
\def\ZZ{\mathop{\bf Z\kern 0pt}\nolimits}
\def\TT{\mathop{\bf T\kern 0pt}\nolimits}
\def\P{\mathop{\bf P\kern 0pt}\nolimits}

\def\St{{\rm S}}

\newcommand{\ds}{\displaystyle}

\newcommand{\Z}{{\mathbf Z}}

\def\dive{\mathop{\rm div}\nolimits}

\newcommand{\andf}{\quad\hbox{and}\quad}
\newcommand{\with}{\quad\hbox{with}\quad}
\newcommand{\beq}{\begin{equation}}
\newcommand{\eeq}{\end{equation}}
\newcommand{\ben}{\begin{eqnarray}}
\newcommand{\een}{\end{eqnarray}}
\newcommand{\beno}{\begin{eqnarray*}}
\newcommand{\eeno}{\end{eqnarray*}}
\newtheorem{defi}{Definition}[section]
\newtheorem{theo}{Theorem}
\newtheorem{lemma}{Lemma}[section]

\newtheorem{corol}{Corollary}[section]
\newtheorem{prop}{Proposition}[section]

\begin{document}

\title[Global solutions to the Navier-Stokes equations in~$\R^{3}$]{Wellposedness and stability results for
  the   Navier-Stokes equations in~$\R^{3}$}

\author[J.-Y. Chemin]{Jean-Yves  Chemin}
\address[J.-Y. Chemin]%
{ Laboratoire J.-L. Lions UMR 7598\\ Universit{\'e} Paris VI\\
175, rue du Chevaleret\\ 75013 Paris\\FRANCE }
\email{chemin@ann.jussieu.fr }
\author[I. Gallagher]{Isabelle Gallagher}
\address[I. Gallagher]%
{ Institut de Math{\'e}matiques de Jussieu UMR 7586\\ Universit{\'e} Paris VII\\
175, rue du Chevaleret\\ 75013 Paris\\FRANCE }
\email{Isabelle.Gallagher@math.jussieu.fr}


\begin{abstract}
In~\cite{cgens}   a class of initial data to the three dimensional, periodic,
incompressible Navier-Stokes  equations was presented, generating a global smooth solution although the norm of
the initial data may be chosen arbitrarily large. The aim of this article is twofold. First, we adapt the construction
of~\cite{cgens} to the case of the whole space: we prove that if a certain nonlinear function of the initial data is small
enough, in a Koch-Tataru~\cite{kochtataru} type space, then there is a global solution to the Navier-Stokes  equations. We
provide an example of initial data satisfying that nonlinear smallness condition, but whose norm is arbitrarily
large in~$  C^{-1}$. Then we prove a stability result on the nonlinear smallness
 assumption. More precisely we show that the new smallness assumption  also holds for
 linear superpositions of
translated and dilated iterates of the initial data, in the spirit of a construction
in~\cite{bcghardy}, thus
generating a large number of different examples.
\end{abstract}


\maketitle

\setcounter{equation}{0}
\section{Introduction}
\subsection{On the global wellposedness of the Navier-Stokes system}
We consider   the three dimensional, incompressible Navier-Stokes system in~$\R^3$,
\[
(NS)\ \left\{
\begin{array}{c}
\partial_{t} u -\Delta u +u\cdot\nabla u=-\nabla p\\
\dive u =0\\
u_{|t=0}=u_{0}.
\end{array}
\right.
\]
Here~$ u$ is a three-component vector field~$ u = (u_{1},u_{2},u_{3})$ representing the
velocity of the fluid, ~$p$ is a scalar denoting the pressure, and both are unknown functions of the space
variable~$ x \in \R^3$ and
of the time variable~$ t \in \R^+$.
We have chosen the kinematic viscosity of the fluid  equal to one
for simplicity -- a comment on the dependence of our results on viscosity is given further down in this introduction.

It is well-known that~$(NS)$ has a global, smooth solution if the initial data is small
enough in the scale invariant space~$\dot
H^{\frac12}$, where we recall that~$\dot H^{s}$ is  the set of tempered distributions~$f$
with Fourier transform~$\widehat f$ in~$L^{1}_{loc}(\R^{3})$ and such that
$$
\|f\|_{\dot H^{s}} \eqdefa \left(
\int_{\R^{3}} |\xi|^{2s} |\widehat f (\xi)|^{2} \: d\xi
\right)^{\frac12}
$$
is finite.
We recall that the scaling of~$(NS)$ is the following:  for any positive~$\lambda$,
the vector field~$u$ is a solution associated with the data~$u_0$ if~$ u_\lambda$
is a solution associated with~$u_{0,\lambda}$, where
$$
 u_\lambda(t,x) =\lambda u(\lambda^2 t,\lambda x) \: \: \: \:
\mbox{and} \: \: \: \: u_{0,\lambda}(x) = \lambda u_0(\lambda
x).
$$
The result in~$\dot
H^{\frac12}$ is due to H. Fujita and T. Kato in~ \cite{fujitakato} (see
also~\cite{leray} for a similar result, where
the smallness of~ $u_{0}$  is measured by~$\|u_{0}\|_{L^{2}}\|\nabla
u_{0}\|_{L^{2}}$). Since then, a number
of works have been devoted to proving similar wellposedness results for larger classes
of initial data; one should mention the result of T. Kato\ccite{kato} where the smallness
is measured in~$L^3$ (see also\ccite{gigamiyakawa}) and the result of M. Cannone, Y.
Meyer and F. Planchon (see\ccite{cannonemeyerplanchon}) where the smallness is
measured in the Besov space~$\dot B^{-1+\frac 3 p}_{p,\infty}$. 
Let us recall that, for positive~$\s$,
\[
\|u\|_{\dot B^{-\s}_{p,r}} \eqdefa 
\Bigl\|t^{\frac {\s} 2} \|\St(t)u \|_{L^{p}}\Bigr\|_{L^{r}(\R^+,\frac {dt} t)}
\]
where~$\St(t) = e^{t\Delta}$
 denotes the heat flow. The importance of this result  can be illustrated by the following
 example: if~$\f$ is a function in the Schwartz space~$\cS(\R^3)$, let us introduce the  family of divergence free vector fields \[ \f_{\e}(x) = \cos \Bigl( \frac {x_{3}}\e\Bigr) (\partial_{2}\f,-\partial_{1}\f,0).
\]
Then, for small~$\e$, the size of~$\|\f_{\e}\|_{\dot B^{-\s}_{p,r}}$ 
is~$\e^{\s}$.

Let us also mention  the result by H. Koch and D.
Tataru in~\cite{kochtataru} where the smallness is measured in the space~$BMO^{-1}$,
defined by
\[
\|u\|_{BMO^{-1}} \eqdefa \|u\|_{\dot B^{-1}_{\infty,\infty}} +
\supetage {x\in \R^3}{R>0}  R^{-\frac 3 2}\Bigl(\int_{P(x,R)}
|\St(t)u (y)|^2 \: dydt\Bigr)^{\frac 1 2},
\]
where~$P(x,R) = [0,R^{2}] \times B(x,R)$ and~$ B(x,R)$ denotes the ball of radius~$R$ centered at zero.

As observed by H. Koch and D. Tataru,
this norm seems to be the ultimate norm for the initial data for  which the classical Picard's
iterative  scheme can work. Indeed the first iterate,~$\St(t)u_{0}$ must
be  in~$L^2$ locally in~$\R^+\times\R^3$. In particular, ~$\St(t)u_{0}$ must be
in~$L^2([0,1]\times B(0,1))$. Then considering the norm of the space must be invariant
by translation as well as by the scaling of the equation, we get the
norm~$\|\cdot\|_{BMO^{-1}}$. Moreover, let us notice that we have
\[
\supetage {x\in \R^3}{R>0}  R^{-\frac 3 2}\Bigl(\int_{P(x,R)}
|\St(t)u(y)|^2dy\Bigr)^{\frac 1 2} \leq \|\St(t) u\|_{L^2(L^\infty)}.
\]
and thus~$\|u\|_{\dot B^{-1}_{\infty,\infty}}\leq \|u\|_{BMO^{-1}}\leq
\|u\|_{\dot B^{-1}_{\infty,2}}$.

Moreover the space~$\dot C^{-1} =  \dot B^{-1 }_{\infty,\infty}$ seems to be optimal independently of
the method of resolution, due to the following   argument (see~\cite{adt} for instance). Let~$B$ be  a Banach space
continuously included in the space~$\cS'$ of tempered distributions on~$\R^{3}$. Let us assume
that, for any~$(\lam,  a)\in \R^+_{\star}\times \R^3$,
\[
\|f(\lam(\cdot-  a))\|_{B}  = \lam^{-1}\|f\|_{B}.
\]
Then we have that~$|\langle f,e^{-|\cdot|^2}\rangle|\leq C\|f\|_{B}$. By
dilation and translation,  we deduce that
\[
\|f\|_{\dot C^{-1}}= \sup_{t>0}t^{\frac 1 2}\|\St(t)f\|_{L^\infty}\leq
C\|f\|_{B}.
\]
We have proved that any Banach space included in~$\cS'$, translation invariant and which has
the right scaling is included in~$\dot C^{-1}$.

Let us  point out that none of the results mentioned so far  are specific to the
Navier-Stokes equations, as they do not use the special structure of the nonlinear
term in~$(NS)$.

Our aim in this paper is to go beyond the smallness condition on the initial data and to
exhibit arbitrarily large initial data in~$\dot C^{-1}  $ which   generate a unique,
global solution.
This was performed in~\cite{cgens} in the periodic case, where we presented a new,
nonlinear smallness assumption
on the initial data, which may hold despite the fact that  the data is   large. That result uses
the structure of the nonlinear term, as it
is based on the fact that the two dimensional Navier-Stokes equation is globally well
posed.

\medbreak
The first theorem of this paper consists in a result of global existence  under a non linear
smallness hypothesis (Theorem\refer{nonlinearcondition} below). The proof consists mainly
in introducing an idea of\ccite{chemin20} in the proof of the Koch and Tataru
Theorem. The non linear smallness hypothesis is, roughly speaking, that the first
iterate~$\St(t)u_{0} \cdot \nabla \St(t)u_{0}$ is exponentially small  with respect
to~$\|u_{0}\|^4_{\dot B^{-1}_{\infty,2}}$.

Then we exhibit an example of a family of initial data with very large~$\dot C^{-1}$
norm which satisfies the non linear smallness hypothesis. This example fits the structure
of the non linear term~$u \cdot \nabla u$.

Then, we study the stability of this nonlinear smallness condition, but not in the usual 
sense of a perturbation by a small vector field. This problem has been solved 
by I. Gallagher, D. Iftimie and F. Planchon in\ccite{gip} and by P. Auscher, 
S. Dubois and P. Tchamitchian in\ccite{adt}. These authors proved that, in any adapted 
scaling space (for instance~$\dot H^{\frac 12 }$,~$L^3$, 
$\dot B^{-1+\frac 3 p}_{p,\infty}$ or~$BMO^{-1}$) the set of initial data giving
rise to global solution is open. 

Our purpose is different. Once constructed an initial data generating a global solution,
we want to generate  a large family of global solutions that may not be   close to the one 
we start with, in the~$\dot C^{-1}$ norm. This is done with a fractal type transform (see the forthcoming 
Definition\refer{defTLamX}). Roughly speaking, this is the linear superposition of 
an arbitrarily large number of dilated and translated 
 iterates of  the initial data, and we will see that the initial data so-transformed
    still satisfies the nonlinear smallness assumption. That of course enables one to construct a
very large class of   initial data
satisfying that smallness assumption; the transformation is based on
 a construction of~\cite{bcghardy}.

\subsection{Definitions} Before  presenting more precisely the results of this paper,
let us give some definitions
and notation.  We shall be using Besov spaces,
which are defined equivalently
using the Littlewood-Paley decomposition or the heat operator. As both definitions will
be useful in the following,
we present them both in the next definition.

\begin{defi}
\label{besov}
{ \sl Let $\varphi \in \mathcal{S}(\R^{3})$ be such that
$\widehat\varphi(\xi) = 1$ for $|\xi|\leq 1$ and $\widehat\varphi(\xi)= 0$ for
$|\xi|>2$. Define, for~$j \in \Z$, the function~$\varphi_{j}(x)\eqdefa
2^{3j}\varphi(2^{j}x)$, and the Littlewood--Paley operators
$ S_{j} \eqdefa \varphi_{j}\ast\cdot \quad \mbox{and} \quad \Delta_{j}
\eqdefa S_{j+1} - S_{j}.
$ Let~$f$ be in~$\mathcal{S}'( \R^{3})$. Then~$f$ belongs to the homogeneous Besov
space~$\dot B^{s}_{p,q}(\R^3)$ if and only if
\begin{itemize}
\item The partial sum $ \sum^{m}_{-m} \Delta_{j}f$ converges towards $f$ as
a tempered distribution;
\item The sequence $\varepsilon_{j} \eqdefa 2^{js}\| \Delta_{j}
f\|_{L^{p}}$ belongs to $\ell^{q}(\Z)$.
\end{itemize}
In that case
$$
\|f\|_{\dot B^{s}_{p,q} } \eqdefa \left(\sum_{j \in \Z}2^{jsq}
\| \Delta_{j}f\|_{L^{p}}^{q}
\right)^{\frac1q}
$$
and if~$s <0$, the one  has the equivalent norm
\beq
\label{besovheatequiv}
\|f\|_{\dot B^{s}_{p,q} } \sim \left\|
t^{-\frac s2} \|\St (t) f\|_{L^{p}}
\right\|_{L^{q}(\R^{+};\frac{dt}t)}.
\eeq}
\end{defi}
Let us notice that the above equivalence comes from the inequality, proved for instance
in\ccite{chemin20},
\beq
\label{Bernsteinexp}
\|\St(t)\Delta_{j} a\|_{L^p} \leq Ce^{-C^{-1}2^{2j}t}
\|\Delta_{j} a\|_{L^p}.
\eeq
Note that the following Sobolev-type continuous embeddings hold:
$$
\dot B^{s_{1}}_{p_{1},q_{1}} \subset \dot B^{s_{2}}_{p_{2},q_{2}}, \quad \mbox{as soon
as} \quad
s_{1} - \frac d{p_{1}} = s_{2} - \frac d{p_{2}} \quad  \mbox{with} \quad  p_{1 }
\leq p_{2}  \quad \mbox{and}\quad  q_{1 }\leq q_{2} .
$$
We shall denote by~$\P$ the Leray projector onto divergence free vector fields
$$
\P = \mbox{Id} - \nabla
\Delta^{-1} \mbox{div} .
$$
Before stating the first result of this paper, let us introduce the following space.
\begin{defi}
\label{espaceBesovKT}
{\sl We shall denote by~$E$ the space of functions~$f$ in~$L^1(\R^+;\dot B^{-1}_{\infty,1})$
such that
\[
\sum_{j\in \ZZ} 2^{-j}\Bigl\|\|\Delta_{j}f(t)\|_{L^\infty}
\Bigr\|_{L^2(\R^+;tdt)}<\infty
\]
equipped with the norm
\[
\|f\|_{E}\eqdefa
\|f\|_{L^1(\R^+;\dot B^{-1}_{\infty,1})} +
\sum_{j\in \ZZ} 2^{-j} \Bigl\|\|\D_{j}f(t)\|_{L^\infty}\Bigr\|_{L^2(\R^+;tdt)}.
\]}
\end{defi}
Let us remark that, for any homogeneous function~$\s$ of order $0$
smooth  outside~$0$, we have
\[
\forall p \in [1,\infty]\, ,\ \|\s(D)\Delta_{j}f \|_{L^p}\leq C \|\Delta_{j}f \|_{L^p}.
\]
Thus the Leray  projection~${\mathbf P}$ onto divergence free vectors fields maps continuously~$E$
into~$E$.

\subsection{Statement of the results}
\subsubsection{Global existence results}
The first result we shall prove is a new global wellposedness result, under a nonlinear smallness assumption on
the initial data.
\begin{theo}\label{nonlinearcondition}
{\sl There is a constant~$C_{0}$ such that the following result holds. Let~$u_{0}
\in \dot H^{\frac12}(\R^{3})$ be a divergence free vector field. Suppose that
\begin{equation}\label{nonlinear}
\Bigl\|\P \Bigl( \St (t)u_{0} \cdot \nabla  \St (t)u_{0}\Bigr)\Bigr\|_E
\leq C_{0}^{-1} \exp \Bigl(
-C_{0} \|u_{0}\|_{\dot B^{-1}_{\infty,2}}^{4} \Bigr).
\end{equation}
Then there is a unique, global solution to~$(NS)$ associated with~ $u_{0}$, satisfying
$$
u \in C_{b}(\R^{+};\dot H^{\frac12}) \cap L^{2}(\R^{+};\dot H^{\frac32}) .
$$}
\end{theo}
{\bf{Remarks} }
\begin{itemize}
\item
As in~\cite{cgens}, Condition~(\ref{nonlinear}) is a nonlinear smallness
condition on the initial data. In particular Theorem~\ref{example}
below provides a class of examples of arbitrarily large vector fields in~$\dot
B^{-1}_{\infty,\infty}$
satisfying~(\ref{nonlinear}).

\item
The proof of Theorem~\ref{nonlinearcondition} is given in Section~\ref{proof1} below; it
consists in writing the solution~$u$ (which exists for a short time at least), as~$u =
\St (t)u_{0} + R$ and in proving
a global wellposedness result for the perturbed Navier-Stokes equation satisfied by~$R$, under
assumption~(\ref{nonlinear}).
\end{itemize}
Now let us give an example of   large initial data satisfying the assumptions of
Theorem~\ref{nonlinearcondition}.
\begin{theo}\label{example}
{\sl Let~ $\phi \in {\mathcal S}(\R^{3})$ be a given function, and consider two real
numbers~$\varepsilon$ and~$\alpha$
in~$]0,1[$. Define
$$
\varphi_{\varepsilon}(x) = \frac{({-\log \varepsilon})^{\frac15}}
{\e^{1-\alpha}} \cos  \left(
\frac{x_{3}}{\e}\right)
\phi\Bigl(x_{1}, \frac{x_{2}}{\varepsilon^{\alpha}}, x_{3}\Bigr).
$$
There is a constant~$C>0$ such that for~$\varepsilon$ small enough, the smooth,
divergence free vector field
$$
u_{0,\varepsilon}(x) = (\partial_{2}\varphi_{\varepsilon}(x),
-\partial_{1}\varphi_{\varepsilon}(x), 0)
$$
satisfies
$$
C^{-1} ({-\log \varepsilon})^{\frac15}
 \leq \|u_{0,\varepsilon}\|_{\dot B^{-1}_{\infty,\infty}}
  \leq C ({-\log \varepsilon})^{\frac15},
$$
and
\beq
\label{Stuepso}
\|\St (t)u_{0,\varepsilon} \cdot \nabla \St (t)u_{0,\varepsilon} \|_{E} \leq C
\varepsilon^{\frac \alpha 3}({-\log \varepsilon})^{\frac25}.
\eeq
Thus for~$\varepsilon$ small enough, the vector field~ $u_{0,\e}$ generates a
unique, global solution to~$(NS)$.
}\end{theo}
The proof of Theorem~\ref{example} is the purpose of Section~\ref{proof2}.

{\bf Remark} One can also write this example in terms of the Reynolds
number of the fluid: let~${\rm Re} > 0$ be the Reynolds
number, and define the rescaled velocity field~$\displaystyle v(t,x ) =\nu
 u (\nu t , x )$ where~ $\nu =  1 / {\rm Re}$.  Then~$ v$ satisfies the
Navier-Stokes equation
$$
\partial_t v + {\bf P} (v\cdot \nabla v) - \nu  \Delta v= 0
$$
and Theorem~\ref{example} states the following: the vector field
$$
 v_{0,\nu} (x)= (- \log \nu)^{\frac15} \cos \left(\frac{x_{3}}\nu\right)
\left(
(\partial_{2} \phi)
(x_{1},\frac{x_{2}}{\nu^{\alpha}}, x_{3}), \nu^{\alpha}
(-\partial_{1} \phi)
(x_{1},\frac{x_{2}}{\nu^{\alpha}}, x_{3}) \right)
$$
satisfies
$$
\| v_{0,\nu}\|_{\dot B^{-1}_{\infty,\infty}} \sim C \nu  (- \log \nu)^{\frac15}
$$
and generates a global solution to the Navier-Stokes equations if~$\nu$ is small enough. Compared with the usual theory of global
existence for the Navier-Stokes equations, we have gained a (power of a) logarithm in the smallness assumption in terms
of the viscosity, since classically
one expects the initial data to be small with respect to~$\nu$.

\subsubsection{Stability results} The second aim of this paper is to give some stability properties of global solutions. It is known since~\cite{gip} that
any initial data in~$\dot B^{-1+\frac3p}_{p,\infty}$ giving rise to a unique global solution is stable: a small perturbation of that data also generates a global solution (see~\cite{adt} for the case of~$BMO^{-1}$).
Here we     present a stability result where the perturbation is as large as the initial data but has a special form: it
consists in the superposition of dilated and translated duplicates of the initial data, in the spirit of
profile decompositions of P. G\'erard (see~\cite{pgprofils}). This
  transform is a version of the fractal transform used in\ccite{bcghardy}
in the study of refined Sobolev and Hardy inequalities.
Let us be more precise and define the transformation. We shall only be considering compactly supported intial
data for this study, and up to a rescaling we shall suppose to simplify that the support of the initial data
is restricted to  the unit cube~$Q$ of~$\R^{3}$ centered at~0.
\begin{defi}\label{defTLamX}
{\sl
Let~$X = (x_{1, \dots,x_{K}}))$ be a set of~$K$ distinct points in~$\R^{3}$.  For~$\Lam \in 2^{\N}$, let us define
\[
\tlk \quad  \left\{
\begin{array}{ccl}
\cS' & \rightarrow & \cS'\\
f & \mapsto & \displaystyle \tlk f\eqdefa \sum_{J\in \{1,\dots , K\}} T_\Lam^{J}f
\end{array}
\right.\with
 T_\Lam^{J}f (x)\eqdefa  \Lam  f(\Lam (x-x_{J}) ).
\]
}
\end{defi}
It can be noted that this is a generalization of the fractal transformation~$T^{k}$   studied in~\cite{bcghardy}.

The next statement is quite easy to prove: it shows that this transformation on the initial data preserves  global wellposedness, as soon as the scaling parameter~$\Lam$ is large enough (the threshold~$\Lambda$ being   unknown as a  function
of the initial data). The theorem following that statement gives a quantitative approach to that stability: if the initial
data~$u_{0} $ satisfies the smallness assumption~(\ref{nonlinear}) of Theorem~\ref{nonlinearcondition}, then so does~$\tlk u_{0}$ as soon
as~$\Lam$ is large enough (the threshold being an explicit function of norms of~$u_{0}$).

More precisely we have the following results.
\begin{prop}\label{stabildebil}
{\sl Let~$u_{0}$ be a divergence free vector field  in~$\dot H^{\frac 1 2}(\R^3)$
generating a  unique, global solution to the Navier-Stokes equations and~$X$ be a
finite sequence of distinct points. There is~$\Lam_{0} >0$ such that, for
any~$\Lam \geq \Lam_{0}$, the vector field~$\tlk u_{0}$ also generates a unique, global
solution.
}
\end{prop}
{\bf Remarks}

\begin{itemize}

\item Using the global stability of global solutions proved in\ccite{gip},  a global solution
associated to an initial  data in~$\dot H^{\frac 1 2}$ is always
in~$L^\infty(\R^+;\dot H^{\frac 1 2})\cap L^2(\R^+;\dot H^{\frac 3 2})$.

\item As the proof of that result in Section~\ref{proofstabildebil}
will show (see Proposition~\ref{stabildebilgeneral}), Proposition~\ref{stabildebil} can be
generalized to the case where   the vector field~$u_{0}$ is
   replaced by any finite sequence   of vector fields in~$\dot H^{\frac12}$ generating a global solution.

\item As we shall see in the proof of Theorem~\ref{stability} stated below, the functions~$\tlk u_{0}$
and~$u_{0}$ have essentially the same norm in~$\dot C^{-1}$.

\end{itemize}
Now let us state the quantitative stability theorem, in particular in the case of an initial data
satisfying the assumptions of Theorem~\ref{nonlinearcondition}. In order to avoid
excessive heaviness, we shall assume from now on that the initial data is compactly supported, and
after scaling, supported in the unit cube~$Q = ]-\frac 12, \frac12 [^{d}$. We shall consider sequences~$X$ such that
\beq\label{distancedelta}
\infetage{(J,J')\in \{1,\dots,K\}2}{J\not =J'} \left\{
d(x_{J},x_{J'} )  ;  d(x_{J},{}^{c}Q)
\right\} \geq \delta>0.
\eeq
We shall prove the following theorem.
\begin{theo}\label{stability}
{\sl
Let~$u_{0}$ be a smooth~$\dot H^{\frac 1 2}$  divergence free vector field,
compactly supported in the cube~$Q$. Suppose that~$u_{0}$ satisfies~(\ref{nonlinear}) in the following slightly looser sense: there is~$\eta \in ]0,1[$ such that
\begin{equation}\label{nonlineareps}
\| {\bf P}(\St(t)   u_{0}\cdot \nabla \St(t) u_{0}) \|_{E} \leq
 C_{0} ^{-1}
\exp \left(
-C_{0}  \bigl(\|u_{0}\|_{\dot B^{-1}_{\infty,2}}+\eta\bigr) ^{4}
\right) - \eta.
\end{equation}
Then there is a positive~$\Lam_{0}$ (depending only on~$ \eta, K,\delta, \|u_{0}\|_{\dot H^{-1}}$
and~$\|u_{0}\|_{\dot B^{-3}_{\infty,\infty}}$)
such that  for
any~$\Lam \geq \Lam_{0}$, the vector field~$   \tlk u_{0}$
 satisfies~(\ref{nonlinear}) and in particular generates a global solution
to~$(NS)$. Moreover, for all~$ r$ in~$[1,\infty]$,
$$
   \|u_{0}\|_{\dot B^{-1}_{\infty,r}} -\eta\leq \|\tlk u_{0 }
\|_{\dot B^{-1}_{\infty,r}} \leq  \|u_{0}\|_{\dot B^{-1}_{\infty,r}} +\eta.
$$}
\end{theo}
\noindent{\bf Remarks }
\begin{itemize}
\item
The factor~$\eta$ appearing in~(\ref{nonlineareps})
means that~$u_{0}$  must not saturate the nonlinear
smallness assumption~(\ref{nonlinear}) of Theorem~\ref{nonlinearcondition}.

\item
The proof  of this  theorem is based on the fact that
the    Besov norm of index~$-1$ as well as~$ \|{\bf P}(\St(t)   u_{0}\cdot \nabla \St(t) u_{0}) \|_{E}$ are
invariant under the action of~$T_{\Lam,X}$, up to some small error terms.
\end{itemize}

As a conclusion of this introduction, let us state the following result, which describes
the action of~$T_{\Lam,X}$ on the family~$u_{0,\e}$ introduced in
Theorem\refer{example}.
\begin{theo}\label{actionfractalonexample}
{\sl Let~$u_{0,\varepsilon}$ be the family introduced in Theorem\refer{example}.
For any~$K$ and~$\delta$, a  constant~$\Lam_{0}$ exists, which is independent of~$\e$, such that the
following result holds. For any family~$X$  and any~$\Lam\geq \Lam_{0} $,    there is  a global solution smooth
solution  of~$(NS)$ with initial  data~$\tlk u_{0,\e}$.
}
\end{theo}

{\bf Remark} Let us point out that as opposed to Proposition~\ref{stabildebil}, Theorem~\ref{stability}
(or rather Lemmas~\ref{technicaltlamk} and~\ref{heatheat} which are the key to its proof)
provides precise bounds on~$\Lam_{0}$ so that the constant~$\Lam_{0}$ appearing in Theorem~\ref{actionfractalonexample}
may be chosen independently of~$\e$.

\section{Proof of Theorem~\ref{nonlinearcondition}} \label{proof1}
\setcounter{equation}{0}
\subsection{Main steps of the proof}
Let us start by remarking that in the case when~$u_{0}$ is small then there is nothing to be proved, so
in the following we shall suppose that~$\|u_{0}\|_{\dot B^{-1}_{\infty,2}}$ is not small,
say~$\|u_{0}\|_{\dot B^{-1}_{\infty,2}} \geq 1$.

We follow the method  introduced by H. Koch and D. Tataru  in\ccite{kochtataru} in order
to look for the solution~$u$ under the form~$ u_{F} + R$, where~$u_{F}(t)
 \eqdefa \St (t)u_{0}$.  Let us denote by~$\cQ$ the bilinear operator defined by
 \[
 \cQ(a,b)(t)\eqdefa -\frac  1 2 \int_{0}^t\St (t-t'){\bf P}\left(a(t')\cdot \nabla b(t')+
 b(t')\cdot\nabla a(t')\right)dt'
 \]
 Then~$R$ is the solution of
 \[
  (MNS)\quad\quad R=\cQ(u_{F},u_{F})+2\cQ(u_{F},R)+\cQ(R,R)  .
 \]
To prove  the global existence of~$u$, we are reduced to proving the global wellposedness of~$(MNS)$; that
  relies on the following easy lemma, the proof of which is omitted.
\begin{lemma}
\label{cacciopoli+}
{\sl Let~$X$ be a Banach space, let~$L$ be a  continuous linear map
from~$X$ to~$X$,  and let~$B$ be a bilinear map from~$X\times X$ to~$X$.
Let us define
\[
\|L\|_{\cL(X)} \eqdefa \sup_{\|x\|=1} \|Lx\|\quad\hbox{and}\quad
\|B\|_{\cB(X)} \eqdefa \sup_{\|x\|=\|y\|=1} \|B(x,y)\|.
\]
If~$\|L\|_{\cL(X)}<1$, then for any~$x_{0}$ in~$X$ such that
\[
\|x_{0}\|_{X}< \frac{ (1-\|L\|_{\cL(X)})^2} {4\|B\|_{\cB(X)} }\virgp
\]
the equation
\[
x=x_{0}+Lx+B(x,x)
\]
has a unique solution in the ball of center~$0$ and
radius~$\ds \frac {1-\|L\|_{\cL(X)}}{2\|B\|_{\cB(X)}}\cdotp$}
\end{lemma}
Let us introduce the functional  space  for which we shall apply the above lemma.
We define the quantity
\[
U(t)\eqdefa \|u_{F}(t)\|_{L^\infty}^2+ t \|u_{F}(t)\|_{L^\infty}^4,
\]
which satisfies
\begin{eqnarray}
\label{prooftheo1eq1}
\int_{0}^\infty U(t)dt & \leq & C \|u_{0}\|_{\dot B^{-1}_{\infty,2}}^2 +
C \|u_{0}\|_{\dot B^{-1}_{\infty,4}}^4 \nonumber \\
 & \leq &C \|u_{0}\|_{\dot B^{-1}_{\infty,2}}^4
\end{eqnarray}
recalling that we have supposed that~$\|u_{0}\|_{\dot B^{-1}_{\infty,2}} \geq 1$ to simplify the notation.

For all~$\lam\geq 0$, let us denote by~$X_{\lam}$ the set of functions on~$\R^+\times \R^3$
such that
\beq
\label{defXlam}
\|v\|_{\lam}\eqdefa \sup_{t > 0}
\biggl( t^{\frac 1 2}\|v_{\lam}(t)\|_{L^\infty}
+\supetage{x\in \R^3}{R>0} R^{-\frac 3 2}
\Big(\int_{P(x,R)}|v_{\lam}(t,y)|^2 dy\Bigr)^{\frac 1 2}\biggr) < \infty ,
\eeq
where
$$
 v_{\lam}(t,x)\eqdefa v(t,x) \exp{\Bigl(- \lam\int_{0}^t  U(t')dt'\Bigr)}
$$
while~$P(x,R)=[0,R^2]\times B(x,R)$
and~$B(x,R)$ denotes the ball  of~$\R^3$ of center~$x$ and radius~$R$.
Let us point out that, in the case when~$\lam=0$, this is exactly  the space introduced by
H. Koch and D. Tataru in\ccite{kochtataru},  and that for any~$\lam\geq 0$  we have due to~(\ref{prooftheo1eq1}),
\beq
\label{prooftheo1eq2}
\|v\|_{\lam}\leq \|v\|_{0} \leq  C\|v\|_{\lam}\exp \Big(
C\lam \|u_{0}\|_{\dot B^{-1}_{\infty,2}}^4
\Bigr).
\eeq
From Lemmas 3.1 and   3.2 of\ccite{kochtataru} together with the above equivalence
of norms, we infer that
\beq
\label{prooftheo1eq3}
\|\cQ(v,w)\|_{\lam} \leq C \|v\|_{\lam}\|w\|_{\lam}
\exp\biggl(C\lam \|u_{0}\|_{\dot B^{-1}_{\infty,2}}^4 \biggr).
\eeq
Theorem\refer{nonlinearcondition} follows from  the following two lemmas.
\begin{lemma}
\label{lemma2prooftheo1}
{\sl There is a constant~$C>0$ such that the following holds. For any non negative~$\lam$, for any~$t \geq 0$
and any~$ f \in E$, we have
\[
\Bigl\|\int_{0}^t \St (t-t')f(t')dt'\Bigr\|_{\lam} \leq C \|f\|_{E}.
\]}
\end{lemma}
\begin{lemma}
\label{lemma3prooftheo1}
{\sl  Let~$u_{0} \in \dot B^{-1}_{\infty,2}$ be given, and define~$u_{F}(t) = \St(t) u_{0}$.
There is a constant~$C>0$ such that the following holds. For any~$\lam\geq 1$,  for any~$t \geq 0$ and
any~$v \in X_{\lam}$, we have
\[
\|\cQ(u_{F},v)(t)\|_{\lam} \leq \frac C {\lam^{\frac 1 4}}\|v\|_{\lam}.
\]}
\end{lemma}
\noindent{\bf End of the proof of Theorem\refer{nonlinearcondition}} Let us apply Lemma~\ref{cacciopoli+}
to Equation~$(MNS)$ satisfied by~$R$, in a space~$X_{\lam}$. We choose~$\lam$ so that according to
 Lemma~\ref{lemma3prooftheo1},
 $$
 \|\cQ(u_{F},\cdot)(t)\|_{{\mathcal L}(X_{\lam})} \leq  \frac14 \cdotp
 $$
 Then according to  Lemma~\ref{cacciopoli+}, there is a unique solution~$R$ to~$(MNS)$ in~$X_{\lam}$
 as soon as~$\cQ(u_{F},u_{F})$ satisfies
 $$
 \|\cQ(u_{F},u_{F})\|_{X_{\lam}} \leq \frac1{16 \|\cQ\|_{{\mathcal B}(X_{\lam})}} \cdotp
 $$
 But~(\ref{prooftheo1eq3}) guarantees that
 $$
 \|\cQ\|_{{\mathcal B}(X_{\lam})}  \leq 
 C \exp\Bigl(C\lam \|u_{0}\|_{\dot B^{-1}_{\infty,2}}^4\Bigr),
 $$
 so it is enough to check that for some constant~$C$,
$$
 \|\cQ(u_{F},u_{F})\|_{X_{\lam}} \leq 
 C^{-1}\exp\Bigl(-C \lam \|u_{0}\|_{\dot B^{-1}_{\infty,2}}^4\Bigr).
 $$
 By Lemma~\ref{lemma2prooftheo1}, this is precisely condition~(\ref{nonlinear}) of
 Theorem\refer{nonlinearcondition}, so under assumption~(\ref{nonlinear}), there is a unique, global solution~$R$
 to~$(MNS)$, in the space~$X_{\lam}$.  This implies immediately that there is a unique, global solution~$u$ to the
 Navier-Stokes system in~$X_{\lam}$.  The fact that~$u$ belongs
 to~$ C_{b}(\R^{+};\dot H^{\frac12}) \cap L^{2}(\R^{+};\dot H^{\frac32})$
 is then simply an argument of propagation of regularity (see for instance~\cite{lemariebook}).

\subsection{ Proof of Lemma\refer{lemma2prooftheo1} }
Thanks to\refeq{prooftheo1eq2}, it is enough to prove Lemma\refer{lemma2prooftheo1} for~$\lam=0$.

Let us start by proving that~$\displaystyle \int_{0}^t \St (t-t')f(t')dt'$ belongs 
to~$L^{2}(\R^+;L^{\infty})$; that
will give in particular the boundedness of the second norm entering in the definition of~$X_{\lam}$.

Using\refeq{Bernsteinexp}, we get
\[
\Bigl\|\int_{0}^t \Delta_{j}\St (t-t')f(t')dt'\Bigr\|_{L^\infty}
\leq C \int_{0}^t  e^{-C^{-1}2^{2j}(t-t')} \|\D_{j}f(t')\|_{L^\infty} dt'.
\]
Young's inequality then gives
\[
\Bigl\|\int_{0}^t \Delta_{j}\St (t-t')f(t')dt'\Bigr\|_{L^2(\R^+;L^\infty)}
\leq C2^{-j} \|\D_{j}f\|_{L^1(\R^+;L^\infty)},
\]
thus the series~$\ds \Bigl( \Delta_{j} \int_{0}^t \St (t-t')f(t')dt'\Bigr)_{j\in\ZZ}$
converges in~$L^2(\R^+;L^\infty)$, and
$$
\label{lemma2prooftheo1eq1}
\Bigl\|\int_{0}^t \St (t-t')f(t')dt'\Bigr\|_{L^2(\R^+;L^\infty)} \leq C \|f\|_{E}.
$$
This implies in particular that
\beq
\label{eq1lemma2prooftheo1}
\supetage{x\in \R^3}{R>0} R^{-\frac 3 2}
\bigg(\int_{P(x,R)}\Bigl|\int_{0}^t (\St (t-t')f(t'))(y)dt'\Bigr|^2 dy\biggr)^{\frac 1 2} \leq C \|f\|_{E} .
\eeq
The second part of the norm defining~$\|\cdot\|_{X_{\lam}}$ in~(\ref{defXlam}) is therefore controlled by the
norm of~$f$ in~$E$.

To estimate the first part of that norm, let us write that for any~$t \geq 0$ and any~$j \in \Z$,
\beno
t^{\frac 1 2} \D_{j}\int_{0}^t \St (t-t')f(t')dt' & = & G_{j}^{(1)}(t)+G_{j}^{(2)}(t)
\quad\hbox{with}\\
G_{j}^{(1)}(t) & \eqdefa & t^{\frac 1 2}\int_{0}^{\frac t 2}  \St (t-t')\Delta_{j}f(t')dt'
\quad\hbox
{and}\\
G_{j}^{(2)}(t) & \eqdefa & t^{\frac 1 2}\int_{\frac t 2}^t  \St (t-t')\Delta_{j}f(t')dt'.
\eeno
 Using again\refeq{Bernsteinexp}  we have, since~$t \leq 2(t-t')$,
 \beno
 \|G_{j}^{(1)}(t) \|_{L^\infty} & \leq & C \int_{0}^{\frac t 2}
 (t-t')^{\frac1 2} 2^je^{-C^{-1}2^{2j}(t-t')} 2^{-j}\|\Delta_{j}f(t')\|_{L^{\infty}}dt'\nonumber\\
 \label{lemma2prooftheo1eq2}  & \leq  &
 2^{-j}\|\Delta_{j}f\|_{L^1(\R^+;L^\infty)}.
 \eeno
 In order to estimate~$\|G_{j}^{(2)}(t)\|_{L^\infty}$,  let us write, since~$t \leq 2t'$,
 \[
 \|G_{j}^{(2)}(t)\|_{L^\infty} \leq C \int_{0}^t e^{-C^{-1}2^{2j}(t-t')}t'^{\frac12}
 \|\D_{j}f(t')\|_{L^\infty}dt'.
 \]
 Using the Cauchy-Schwarz inequality, we get
 \[
 \|G_{j}^{(2)}(t)\|_{L^\infty} \leq C2^{-j}
 \|t^{\frac1 2} \D_{j}f(t)\|_{L^2(\R^+;L^\infty)}.
 \]
  Then using\refeq{eq1lemma2prooftheo1} and summing over~$j \in \Z$ concludes
  the proof of Lemma\refer{lemma2prooftheo1}. \hfill $\square$
\subsection{ Proof of Lemma\refer{lemma3prooftheo1} } We have (see for
instance\ccite{kochtataru} or\ccite{chemin22}) that
\begin{eqnarray*}
\cQ(v,w)(t,x) & = & \int_{0}^t\int_{\R^3} k(t-t',y)v(t',x-y)w(t',x-y)dy dt'\\
 & = & k \star (vw)(t,x) \quad \mbox{with} \quad
|k(\tau,\zeta)| \leq \frac C {(\sqrt \tau +|\zeta|)^4}\cdotp
\end{eqnarray*}
The proof relies  now mainly on the following proposition.
\begin{prop}
\label{propprooflemma2prooftheo1}
{\sl Let~$u_{0} \in \dot B^{-1}_{\infty,2}$ be given, and define~$u_{F}(t) = \St(t)u_{0}$.
There is a constant~$C$ such that the following holds. Consider, for any positive~$R$ and for~$(\tau,\zeta)
\in \R^+\times \R^{3}$, the following functions:
\[
K_{R}^{(1)}(\tau,\zeta)\eqdefa {\bf 1}_{|\zeta|\geq R}\frac 1 {|\zeta|^4}\andf
K_{R}^{(2)}(\tau,\zeta)\eqdefa
{\bf 1}_{|\zeta|\leq R}\frac 1 { (\sqrt \tau+|\zeta|)^4}\cdotp
\]
Then for any~$\lam\geq 1$ and any~$R>0$,
\beq
\label{propprooflemma2prooftheo1eq1}
\Bigl\| e ^{-\lam\int_{0}^tU(t')dt'}K_{R}^{(1)}\star (u_{F}v)
\Bigr\|_{L^\infty([0,R^2]\times\R^3)}
\leq  \frac C { \lam^{\frac 1 2} R }\|v\|_{\lam}.
\eeq
Moreover, for any~$\lam\geq 1$ and  any~$R>0$, \beq
\label{propprooflemma2prooftheo1eq2}
\Bigl\| e ^{-\lam\int_{0}^tU(t')dt'}K_{R}^{(2)}\star (u_{F}v) 
\Bigr\|_{L^\infty([R^2,2R^2] \times \R^3)} \leq
\frac C { \lam^{\frac 1 4} R }\|v\|_{\lam}.
\eeq}
\end{prop}

 \medbreak\noindent{\bf Proof of Proposition\refer{propprooflemma2prooftheo1} }
Let us write that
\beno
V_{\lam}^{(1)}(t,x) & \eqdefa &
e ^{-\lam\int_{0}^tU(t')dt'}|K_{R}^{(1)}\star (u_{F}v)(t,x)| \\
 & \leq  & \int_{0}^t\int_{{}^c B(0,R)}
\frac 1 {|y|^4}e ^{-\lam\int_{t'}^tU(t'')dt''} \|u_{F}(t',\cdot)\|_{L^\infty}
|v_{\lam}(t',x-y)|dt'dy.
\eeno
By the Cauchy-Schwarz inequality and by definition of~$U$, we infer that
\ben
V_{\lam}^{(1)}(t,x) & \leq  & \biggl(\int_{0}^t\int_{{}^c B(0,R)}
\frac 1 {|y|^4}e ^{-2\lam\int_{t'}^tU(t'')dt''} \|u_{F}(t',\cdot)\|^2_{L^\infty}dt'dy
\biggr)^{\frac 1 2}\nonumber\\
 & & \qquad\quad\qquad\quad\qquad\quad{}\times
\biggl(\int_{0}^t\int_{{}^c B(0,R)}
\frac 1 {|y|^4}
|v_{\lam}(t',x-y)|^2 dt'dy \biggr)^{\frac 1 2}\nonumber\\
\label{propprooflemma2prooftheo1eq11}
& \leq  & \biggl(\frac C {\lam R}\biggr)^{\frac 1 2}
\biggl(\int_{0}^t\int_{{}^c B(0,R)}
\frac 1 {|y|^4} |v_{\lam}(t',x-y)|^2 dt'dy\biggr)^{\frac 1 2}.
\een
Now let us decompose the integral on the right on rings; this gives
\beno
\int_{0}^t\int_{{}^c B(0,R)}
\frac 1 {|y|^4} |v_{\lam}(t',x-y)|^2 dt'dy
& = & \sum_{p=0}^\infty \int_{0}^t\int_{B(0,2^{p+1}R)\setminus B(0,2^pR)}
\frac 1 {|y|^4} |v_{\lam}(t',x-y)|^2 dt'dy \\
& \leq  &\frac 1 R\sum_{p=0}^\infty 2^{-p+3} (2^{p+1}R)^{-3} \\
&& \qquad{}\times
\int_{0}^{t}\int_{B(0,2^{p+1}R)}
 |v_{\lam}(t ,x-y)|^2 dt dy .
 \eeno
 As~$t\leq R^2$ and~$p$ is non negative, we have
 \beno
 \int_{0}^t\int_{{}^c B(0,R)}
\frac 1 {|y|^4} |v_{\lam}(t',x-y)|^2 dt'dy
& \leq  & \frac C R\sum_{p=0}^\infty 2^{- p} (2^{p+1}R)^{-3}
\int_{P(x,2^{p+1}R)} |v_{\lam}(t,z)|^2 dtdz\\
& \leq & \frac C R\sum_{p=0}^\infty 2^{- p} \sup_{R'>0} \frac 1{R'^{3}}
\int_{P(x,R')} |v_{\lam}(t,z)|^2 dtdz.
\eeno
By definition of~$\|\cdot\|_{\lam}$, we infer that
$$
\int_{0}^t\int_{{}^c B(0,R)}
\frac 1 {|y|^4} |v_{\lam}(t',x-y)|^2 dt'dy \leq    \frac C R \|v\|_{\lam}^2,
$$
Then, using\refeq{propprooflemma2prooftheo1eq11}, we conclude  the proof 
of\refeq{propprooflemma2prooftheo1eq1}.

In order to prove the second inequality,
 let us observe  that
 \beno
  &&e ^{-\lam\int_{0}^tU(t')dt'}|(K_{R}^{(2)}\star (u_{F}v))(t,x) |   
  \leq  \cK_{R}^{(21)}(t,x)+\cK_{R}^{(21)}(t,x)\with\\
 \cK_{R}^{(22)}(t,x) & \eqdefa  & \int_{0}^{\frac t 2}\int_{B(0,R)}
 \frac {1} {(\sqrt {t-t'}+|y|)^4} e ^{-\lam\int_{t'}^tU(t'')dt''}
 \|u_{F}(t',\cdot)\|_{L^\infty} |v_{\lam}(t',x-y)|  dt'dy\,,\\
  \cK_{R}^{(22)}(t,x) & \eqdefa  & \int_{\frac t 2}^t\int_{B(0,R)}
 \frac {1} {(\sqrt {t-t'}+|y|)^4} e ^{-\lam\int_{t'}^tU(t'')dt''}
 \|u_{F}(t',\cdot)\|_{L^\infty} |v_{\lam}(t',x-y)|  dt'dy.
 \eeno
 Using the Cauchy-Schwarz inequality, as~$t\in [R^2,2R^2]$ and~$t \leq 2(t-t')$, we infer that
 \beno
  \cK_{R}^{(21)}(t,x) &\leq & \biggl(\int_{0}^{\frac t2}\int_{B(0,R)}
   \frac {1} {(\sqrt {t-t'}+|y|)^8} e ^{-2\lam\int_{t'}^tU(t'')dt''}
 \|u_{F}(t',\cdot)\|^2_{L^\infty}dt'dy\biggr)^{\frac 1 2}\nonumber\\
 && \qquad\qquad\qquad\qquad\qquad\qquad\qquad{}\times
 \biggl( \int_{0}^{\frac t 2}\int_{B(0,R)} |v_{\lam}(t',x-y)|^2 dt'dy\biggr)^{\frac 1 2}
 \nonumber\\
  & \leq &
  \frac{C}{\lam^{\frac12}} \left(
  \int_{B(0,R)} \frac{dy}{(R+|y|)^{8}}
  \right)^{\frac12} \biggl( \int_{0}^{\frac t 2}\int_{B(0,R)}
   |v_{\lam}(t',x-y)|^2 dt'dy\biggr)^{\frac 1 2}
   \nonumber\\
 & \leq & \frac C {(t\lam)^{\frac 1 2}}R^{-\frac 3 2}
  \biggl( \int_{0}^{R^2}\int_{B(0,R)} |v_{\lam}(t',x-y)|^2 dt'dy\biggr)^{\frac 1 2},
 \eeno
 so that
 \beq  \label{propprooflemma2prooftheo1eq12}
 K_{R}^{(21)}(t,x)
    \leq     \frac C {(t\lam)^{\frac 1 2}}\|v\|_{\lam}.
\eeq
In order to estimate~$\cK_{R}^{(22)}$, let us write that
\beno
\cK_{R}^{(22)}(t,x) & \leq & \int_{\frac t 2}^t\int_{\R^3}
 \frac {1} {(\sqrt {t-t'}+|y|)^4} e ^{-\lam\int_{t'}^tU(t'')dt''}
 \|u_{F}(t',\cdot)\|_{L^\infty}\|v_{\lam}(t',\cdot)\|_{L^\infty} dt'dy\\
  & \leq  & C\|v\|_{\lam}\int_{\frac t 2}^t
 \frac {1} {\sqrt {t-t'}} e ^{-\lam\int_{t'}^tU(t'')dt''}
 \frac{ \|u_{F}(t',\cdot)\|_{L^\infty}}{t'^{\frac 1 2}} dt' .
\eeno
By definition of~$U$ and using the fact that~$t \leq 2t'$, H\"older's inequality 
implies that
\beno
\cK_{R}^{(22)}(t,x) & \leq  & \frac C {t^{\frac 1 2}}\|v\|_{\lam}
 \Bigl(\int_{0}^t  e ^{-4\lam\int_{t'}^tU(t'')dt''} t'\|u_{F}(t',\cdot)\|^4_{L^\infty}dt'\Bigr)^{\frac 14}\\
& \leq &\frac  C{\lam^{\frac 1 4}t^\frac 1 2}\|v\|_{\lam}.
\eeno
Together with\refeq{propprooflemma2prooftheo1eq12}, this concludes the proof of the
proposition. \hfill $\square$

 From this proposition, we infer immediately the following corollary. This corollary proves
 directly one half of Lemma~\ref{lemma3prooftheo1}, as it gives a control of~$\cQ(u_{F},v)$ in
 the first norm out of the two entering in the definition  of~$X_{\lam}$.
 \begin{corol}
 \label{corprooftheo1}
 {\sl
 Under the assumptions of Proposition~\ref{propprooflemma2prooftheo1}, we have
\[
 t^{\frac12} e ^{-\lam\int_{0}^tU(t')dt'} \|\cQ(u_{F},v)(t,\cdot)\|_{L^\infty}\leq
 \frac C {\lam^{\frac 1 4}}\|v\|_{\lam}.
  \]
 }\end{corol}

\medbreak\noindent{\bf Proof of Corollary\refer {corprooftheo1}}.
 Let us write that
 \[
 k\star( u_{F}v)\,(t,x) = k\star( u_{F}{\bf 1}_{{}^c B(x,2\sqrt t)}v)\,(t,x)
 +k\star (u_{F}{\bf 1}_{B(x,2\sqrt t)}v)\,(t,x).
 \]
From Proposition\refer{propprooflemma2prooftheo1}, we infer that
\beno
e ^{-\lam\int_{0}^tU(t')dt'}|k\star( u_{F}{\bf 1}_{{}^c B(x,2\sqrt t)}v)\,(t,x)| & \leq &
e ^{-\lam\int_{0}^tU(t')dt'} K_{2\sqrt t}^{(1)}\star( |u_{F}{\bf 1}_{B(x,2\sqrt t
)}v|)
(t,x)\\
& \leq & \frac C{( t\lam)^\frac 1 2}\|v\|_{\lam}.
\eeno
Moreover, thanks to Proposition\refer{propprooflemma2prooftheo1}, we have also
\beno
e ^{-\lam\int_{0}^tU(t')dt'}|k\star (u_{F}{\bf 1}_{B(x,2\sqrt t )}v)\,(t,x)| & \leq  &e ^{-\lam\int_{0}^tU(t')dt'}
K_{2\sqrt t} ^{(2)}\star  (|u_{F}|{\bf 1}_{B(x,2\sqrt t )}|v|) \,(t,x)\\
 & \leq & \frac C{\lam^{\frac 1 4} t^{\frac 12}} \|v\|_{\lam}.
\eeno
This proves the corollary.\hfill $\square$

\medbreak In order to conclude the proof of Lemma\refer{lemma3prooftheo1}, let us
estimate~$\|k\star (u_{F}v)\|_{L^2(P(x,R))}$, for an arbitrary~ $x \in \R^{3}$. Let us write that
\[
k\star (u_{F}v)=  k\star (u_{F}{\bf 1}_{{{}^c B}(x,2R)}v)+k\star (u_{F}{\bf 1}_{{B}(x,2R)}v).
\]
Observing that, for any~$y\in B(x,R)$, we have
\[
|k\star (u_{F}{\bf 1}_{{{}^c B}(x,2R)}v) (t,y)|\leq  CK_{R}^{(1)}\star
(|u_{F}|{\bf 1}_{{{}^c B}(x,2R)}|v|) (t,y),
\]
and using Inequality\refeq{propprooflemma2prooftheo1eq1} of
Proposition\refer{propprooflemma2prooftheo1}, we get
\[
\|e ^{-\lam\int_{0}^tU(t')dt'}
k\star (u_{F}{\bf 1}_{{{}^c B}(x,2R)}v)\|_{L^\infty(P(x,R))} 
\leq\frac C {\lam^{\frac 1 2}R}\|v\|_{\lam}.
\]
As the volume of $P(x,R)$ is proportional to~$R^5$, we infer that
\beq
\label {prooftheo1eq14}
\|e ^{-\lam\int_{0}^tU(t')dt'}
k\star (u_{F}{\bf 1}_{{{}^c B}(x,2R)}v)\|_{L^\infty(P(x,R))}  \leq
\frac C {\lam^{\frac 1 2}}R^{\frac 3 2}\|v\|_{\lam}.
\eeq
The following inequality is easy and classical, so its proof is omitted.
\beq
\label {prooftheo1eq15}
\biggl\|e ^{-\lam\int_{0}^tU(t')dt'}\cQ(u_{F},v)(t)\biggr\|_{L^2([0,T]\times\R^3)}\leq
\frac C{\lam^{\frac 1 2}} \|v_{\lam}\|_{L^2([0,T]\times \R^3)}.
\eeq
We deduce that
\beno
\biggl\|e ^{-\lam\int_{0}^tU(t')dt'} k\star (u_{F}{\bf 1}_{B(x,2R)}) v
\biggr\|_{L^2(P(x,R))}
& \leq  & \biggl\|e ^{-\lam\int_{0}^tU(t')dt'} k\star (u_{F}{\bf 1}_{B(x,2R)}v)
\biggr\|_{L^2([0,R^2]\times \R^3)}\\
 & \leq  &   \frac  C{\lam^{\frac 1 2}}
 \|{\bf 1}_{B(x,2R)}v_{\lam}\|_{L^2([0,R^2]\times \R^3)}\\
  &\leq &  \frac  C{\lam^{\frac 1 2}}
 \|v_{\lam}\|_{L^2(P(x,2R))}.
\eeno
 This concludes the proof of Lemma\refer{lemma3prooftheo1}.\hfill $\square$


\section{Proof of Theorem~\ref{example}} \label{proof2}
\setcounter{equation}{0}
In this paragraph we shall check that the vector field~$u_{0,\varepsilon}$ introduced in
Theorem~\ref{example} satisfies the nonlinear smallness assumption of
Theorem~\ref{nonlinearcondition}, and we shall also show that
its~$\dot B^{-1}_{\infty,\infty}$ norm is equivalent
to~$({-\log \varepsilon})^{\frac15}$.  Let us start by proving the following lemma.
\begin{lemma}
\label{fepsilon}
{\sl Let~$f \in {\mathcal S}(\R^{3})$ be given
 and~$\s\in \ds \biggl]0,3\Bigl(1-\frac 1 p\Bigr)\biggr[$. There is
a constant~$C>0$ such that for any~ $\varepsilon \in  ]0,1[$,
the function
$$
f_{\varepsilon}(x) \eqdefa e^{i\frac{x_{3}}\varepsilon}
f\Bigl(x_{1},\frac{x_{2}}{\varepsilon^{\alpha}},x_{3}\Bigr)
$$
satisfies, for  all~$p \geq 1$,
$$
\|f_{\varepsilon}\|_{\dot B^{-\sigma}_{p,1}} \leq C \e^{\sigma + \frac\alpha p}
\andf \|f_{\e}\|_{\dot B^{-\s}_{\infty,\infty}}\geq C^{-1}\e^{\s}.
$$}
\end{lemma}
\begin{proof}
Let us recall that
$$
\|f_{\varepsilon}\|_{\dot B^{-\sigma}_{p,1}} = \sum_{j \in \Z} 2^{-j\sigma}
\|\Delta_{j} f_{\e}\|_{L^{p}}.
$$
We shall  start by estimating the high frequencies, defining a threshold~ $j_{0} \geq 0$ 
to be determined later on. We have
\ben \label{highfreqfeps}
\sum_{j \geq j_{0}} 2^{-j\sigma}
\|\Delta_{j} f_{\e}\|_{L^{p}} & \leq & C2^{-j_{0} \sigma} \| f_{\e}\|_{L^{p}} \nonumber \\
 & \leq & C2^{-j_{0} \sigma} \e^{\frac\alpha p} \|f\|_{L^{p}}.
\een
On the other hand, we have
\beno
\Delta_{j} f_{\e} (x)& =& 2^{3j} \int_{\R^{3}} h(2^{j}(x-y)) f_{\e}(y) \: dy \\
& = & 2^{3j} \int_{\R^{3}} h(2^{j}(x-y))e^{i\frac{y_{3}}\e}
f (y_{1},\frac{y_{2}}{\e^\al},y_{3}) \: dy,
\eeno
so noticing that~$e^{i\frac{y_{3}}\e} = (-i\e \partial_{3})^{N} (e^{i\frac{y_{3}}\e} )$, we get for
any~$N \in \N$,
$$
\Delta_{j} f_{\e} (x) = (i\e)^{N} 2^{3j} 
\sum_{\ell = 0}^{N} C_{N}^{\ell} \int_{\R^{3}} e^{i\frac{y_{3}}\e}
\partial_{3}^{\ell} \bigl(h(2^{j}(x-y))\bigr) \partial_{3}^{N-\ell}
f(y_{1},\frac{y_{2}}{\e^\al},y_{3}) \: dy.
$$
Young's inequality enables us to infer that
$$
2^{-j\sigma} \|\Delta_{j} f_{\e}\|_{L^{p}} \leq C \e^{N} 2^{ j(3-\sigma)} \min\Bigl(
\sum_{\ell = 0}^{N} 2^{j(\ell-3)} \e^{\frac\alpha p}  , \sum_{\ell = 0}^{N}
2^{j(\ell-\frac3p) }\e^{\alpha} \Bigr).
$$
So, choosing~$N$ large enough and since~$\sigma < 3(1-\frac1p)$, we get
\ben \label{lowfreqfeps}
\sum_{j \leq j_{0}} 2^{-j\sigma}
\|\Delta_{j} f_{\e}\|_{L^{p}} & \leq & \sum_{j \leq  {0}} 2^{-j\sigma}
\|\Delta_{j} f_{\e}\|_{L^{p}} +
\sum_{0<j \leq {0}} 2^{-j\sigma}
\|\Delta_{j} f_{\e}\|_{L^{p}} \nonumber \\
& \leq &  C\sum_{ j< 0} 2^{-j(\sigma - 3(1-\frac1p))} \e^{N+\alpha}
+ C\sum_{0 \leq j \leq j_{0}} 2^{j(N-\sigma)} \e^{N+\frac\alpha p} \nonumber \\
& \leq &C \e^{N+\alpha} + C 2^{j_{0}(N-\sigma)} \e^{N+\frac\alpha p}.
\een
Finally choosing~$2^{-j_{0}}= \e $ in~ (\ref{highfreqfeps}) and~(\ref{lowfreqfeps})
ends the proof of the bound
on~$\|f_{\e}\|_{\dot B^{-\s}_{p,1}}$.

In order to go from below~$\|f_{\e}\|_{\dot B^{-\s}_{\infty,\infty}}$,
let us first  observe that, as the space of smooth compactly supported functions is dense
in~$\cS$ and the Fourier transform is continuous on~$\cS$, for any positive~$\eta$, a
function~$g$ exists, the Fourier transform of which is smooth and compactly supported such
that, denoting as before~$\displaystyle g_{\e} (x) = e^{i\frac{x_{3}}\varepsilon}
g(x_{1},\frac{x_{2}}{\varepsilon^{\alpha}},x_{3}
) $,
\beq
\label{prooflemmafepsiloneq1}
\|f_{\e}-g_{\e}\|_{\dot B^{-\s}_{\infty,\infty}}\leq \eta \e^{\s}\andf
\|f-g\|_{L^\infty}\leq \eta.
\eeq
As the support of the Fourier transform of~$g$ is included in the ball~$B(0,R)$ for some
positive~$R$, that
of~$g(x_{1},\e^{-\alpha}x_{2},x_{3})$ is included in the ball~$B(0, R\e^{-\al})$.
Then the support of~$\cF g_{\e}$ is included in the
ball~$B(\e^{-1}(0,0,1), \e^{-\alpha}R)$. This ball is included in~$\e^{-1}\cC$
for some ring~$\cC$. Thanks to\refeq{besovheatequiv} we shall use the heat flow.
Let us write that
\beno
\|g_{\e}\|_{\dot B^{-\s}_{\infty,\infty}} & \sim & \sup_{t>0} t^{\frac \s 2}
\|\St(t)g_{\e}\|_{L^\infty}\\
& \geq & C\e^{\s}\| \St(\e^{2})g_{\e}\|_{L^\infty}.
\eeno
For any function~$h$ such that the support of~$\wh h$ is included in~$\e^{-1}\cC$,
we have
\[
\|\cF^{-1} (e^{\e^2|\xi|^2}\wh h)\|_{L^\infty} \leq C \|h\|_{L^\infty}.
\]
Applied with~$h=\St(\e^2)g_{\e}$, this inequality gives
\[
\|g_{\e}\|_{L^\infty}\leq C \|\St(\e^2)g_{\e}\|_{L^\infty}
\quad\hbox{and thus}\quad
\|g_{\e}\|_{\dot B^{-\s}_{\infty,\infty}}\geq C^{-1}\e^{\s}
\|g_{\e}\|_{L^\infty}=C^{-1}\e^{\s}\|g\|_{L^\infty}.
\]
Now let us write that
\beno
\|f_{\e}\|_{\dot B^{-\s}_{\infty,\infty}} & \geq &
\|g_{\e}\|_{\dot B^{-\s}_{\infty,\infty}} -\eta \e^\s\\
& \geq & C^{-1}\e^\s (\|f\|_{L^\infty}-2\eta).
\eeno
This ends the proof of the lemma.
\end{proof}
This enables us to infer immediately the following corollary.
\begin{corol}
\label{cordonneeinitiale}
{\sl A constant~$C$ exists  such that, for any~$p \geq 3/2$, we have
$$
\|u_{0,\e}\|_{\dot B^{-1}_{p,1}} \leq C
\e^{\frac\alpha p} ({-\log \e})^{\frac15}\quad \mbox{and} \quad
\|u_{0,\e}\|_{\dot B^{-1}_{\infty,\infty}}\geq
C^{-1}({-\log \e})^{\frac15}.
$$}
\end{corol}

The last verification to be made is  the nonlinear assumption~(\ref{Stuepso}). It is
based on the
following lemma.
\begin{lemma}\label{fepsgeps}
{\sl There is a constant~$C$ such that the following result holds.
Let~ $f$ and~$g$ be in~$\dot B^{-1}_{\infty,2} \cap \dot H^{-1}$. Then we have
 $$
\|{\bf P}(\St (t) f\St (t) g)\|_{E} \leq
  C\Bigl(\|f\|_{\dot B^{-1}_{\infty,2}}
   \|g\|_{\dot B^{-1}_{\infty,2}}\Bigr)^{\frac 2 3}
\Bigl ( \|f\|_{\dot H^{-1}} \|g\|_{\dot H^{-1}}\Bigr)^{\frac 1 3}
  $$}
\end{lemma}
\begin{proof}
As the Leray projection~${\bf P}$ is continuous on~ $E$, it is enough to prove the lemma
without~${\bf P}$. Using Bernstein's estimate, we get that
\[
\|\Delta_{j}(\St (t) f \St (t) g)\|_{L^\infty}
\leq C 2^{3j}\|\St (t) f \St (t) g\|_{L^1}.
\]
Then, using the Cauchy-Schwarz inequality, we infer that
\beno
E_{j} & \eqdefa &
\|\Delta_{j}(\St (t) f \St (t) g)\|_{L^1(\R^+;L^\infty)}+
\|t^{\frac 1 2}\Delta_{j}(\St (t) f \St (t) g)\|_{L^2(\R^+;L^\infty)}\\
 & \leq  & C2^{3j}\Bigl(\|\St (t) f \|_{L^2(\R^+;L^2)}
 +\|t^{\frac 1 2}\St (t) f \|_{L^\infty(\R^+;L^2)}\Bigr)
 \|\St (t) g \|_{L^2(\R^+;L^2)}.
\eeno
So using\refeq{besovheatequiv}, we deduce that
\beq
\label{prooffepsgepseq1}
E_{j} \leq  C2^{3j} \|f\|_{\dot H^{-1}}
\|g\|_{\dot H^{-1}}.
\eeq
Let us observe that we also have
\beno
E_{j}
 & \leq  & C\Bigl(\|\St (t) f \|_{L^2(\R^+;L^\infty)}
 +\|t^{\frac 1 2}\St (t) f \|_{L^\infty(\R^+;L^\infty})\Bigr)
 \|\St (t) g \|_{L^2(\R^+;L^\infty)}\\
& \leq & C \|f\|_{\dot B^{-1}_{\infty,2}}
\|g\|_{\dot B^{-1}_{\infty,2}}.
\eeno
Using this estimate for high frequencies and\refeq{prooffepsgepseq1}
for low frequencies, we get, for any~$j_{0}$ in~$\ZZ$,
\beno
\|\St (t)f\St (t)g\|_{E} & =  & \sum_{j}2^{-j}E_{j}\\
 & \leq  &C \biggl(\|f\|_{\dot H^{-1}}\|g\|_{\dot H^{-1}}
 \sum_{j\leq j_{0}}2^{2j}
 +\|f\|_{\dot B^{-1}_{\infty,2}}\|g\|_{\dot B^{-1}_{\infty,2}}
 \sum_{ j\geq j_{0}} 2^{-j}\biggr)\\
& \leq  &
C \Bigl(\|f\|_{\dot H^{-1}}\|g\|_{\dot H^{-1}} 2^{2j_{0}}
 +\|f\|_{\dot B^{-1}_{\infty,2}}\|g\|_{\dot B^{-1}_{\infty,2}}2^{-j_{0}}\Bigr).
\eeno
Choosing~$j_{0}$ such that
\[
2^{3j_{0}}\sim
\frac {\|f\|_{\dot B^{-1}_{\infty,2}}\|g\|_{\dot B^{-1}_{\infty,2}}}
{\|f\|_{\dot H^{-1}}\|g\|_{\dot H^{-1}}
}
\]
 gives the result.
\end{proof}
Finally we are ready to prove estimate~(\ref{Stuepso}). Note that the proof relies
heavily on the special structure
of the nonlinear term in the system. We indeed start by remarking that there is no
derivative in the third direction
since~$u_{0,\e}$
does not have a third component. Then denoting~$u_{F} (t) = \St(t) u_{0,\e}$, we have by
an easy computation and
 with the  notation  as in Lemma\refer{fepsilon},
$$
u_{F}^{1} \partial_{1} u_{F}^{1} + u_{F}^{2} \partial_{2} u_{F}^{1} =
 \frac1{\e^{2}} (-\log \e)^{\frac25}
\St(t)f_{\e}
\St(t)g_{\e}
\quad \mbox{and}$$$$
u_{F}^{1} \partial_{1} u_{F}^{2} + u_{F}^{2} \partial_{2} u_{F}^{2}
= \frac1{\e^{2-\alpha}} (-\log \e)^{\frac25} \St(t)\wt f_{\e}
\St(t)\wt g_{\e},
$$
where~$f$, $\wt f$, $g$ and~$\wt g$ are smooth functions.
The result follows immediately using Lemmas~\ref{fepsgeps}
and Corollary\refer{cordonneeinitiale} together with the fact that the Leray
projection onto divergence free vector fields maps continuously~$E$ into~$E$. \hfill $\square$


\section{Stability results}\label{proof3}
\setcounter{equation}{0}
In this section we shall prove Proposition~\ref{stabildebil}, as well as
 Theorems~\ref{stability} and~\ref{actionfractalonexample} stated in the
introduction. The proof of Proposition~\ref{stabildebil} is rather easy and is given for the sake of completeness
in the next section. The proof of Theorem~\ref{stability} is   the object of Section~\ref{plusdur}
below. Finally Theorem~\ref{actionfractalonexample} is an easy consequence of the methods developped in the
proof of  Theorem~\ref{stability} and is postponed to the end of Section~\ref{plusdur}.

\subsection{Proof of Proposition~\ref{stabildebil}} \label{proofstabildebil}
Proposition~\ref{stabildebil} is an immediate consequence of the following more general result.
\begin{prop}\label{stabildebilgeneral}
{\sl
Let~$X = (x_{1},\dots,x_{K})$ be
  a family of~$K$ distinct   points, and~$(u_{0,1}, \dots, u_{0,K})  $   a family of
divergence free vector fields in~$\dot H^{\frac12}$, each
  generating a unique, global solution to the Navier-Stokes equations.
Then there is~$\Lam_{0} >0$ such that for any~ $\Lam \geq \Lam_{0}$, the vector field
$$
u_{0,\Lam}  \eqdefa  \sum_{J \in \{1,\dots,K\}} T_{\Lam}^{J} (u_{0,J})
$$
also
generates a unique, global solution to the Navier-Stokes equations.
}
\end{prop}
\begin{proof}
The proof of that result is similar to methods of~\cite{bahourigerard} concerning profile decompositions
 (see~\cite{igprofils} for the case of the Navier-Stokes equations). Let us denote by~$u_{J}$ the solution of~$(NS)$
associated with~$u_{0,J}$, and define
$$
u_{\Lam,J} (t,x) = \Lam u_{J} \bigl(\Lam^{2}t, \Lam (x-x_{J})\bigr),
$$
which solves~$(NS)$ with data~$u_{0,\Lam,J} = T_{\Lam}^{J} (u_{0,J}) $.
Then we define the solution~$u_{\Lam}$ of~$(NS)$ with data~$u_{0,\Lam}$, which a priori exists only for a short time. We
 can decompose
 $$
 u_{\Lam} = \sum_{J \in \{1,\dots,K\}} u_{\Lam,J} + R_{\Lam} = u_{\Lam}^{(1)} + R_{\Lam},
 $$
 and~$ R_{\Lam}$ solves the following perturbed Navier-Stokes equation
 $$
 \partial_{t} R_{\Lam} - \Delta R_{\Lam}+ {\mathbf P}(R_{\Lam} \cdot \nabla R_{\Lam}) +
 {\mathbf P}(u_{\Lam}^{(1)} \cdot \nabla R_{\Lam}) + {\mathbf P}(R_{\Lam} \cdot \nabla u_{\Lam}^{(1)} )
 = F_{\Lam}
 $$
 with initial data zero, and where
 $$
 F_{\Lam} = - {\mathbf P} \sum_{J \neq J'} u_{\Lam,J}  \cdot \nabla u_{\Lam,J'}.
 $$
 It is not difficult to prove (see for instance~\cite{igprofils}, Proposition A.2) that~$ R_{\Lam}$ is
 globally defined and unique in~$L^{\infty}(\R^{+};\dot H^{\frac12})
 \cap L^{2}(\R^{+};\dot H^{\frac32}) $ under the condition that
\beq\label{conditionfacile}
 \| F_{\Lam}\|_{L^{2}(\R^{+};\dot H^{-\frac12})} \leq C_{0}^{-1} \exp \left(
 -C_{0} \|u_{\Lam}^{(1)}\|^{4}_{L^{4}(\R^{+};\dot H^{1})}
 \right),
 \eeq
 so let us compute~$ \| F_{\Lam}\|_{L^{2}(\R^{+};\dot H^{-\frac12})} $
  and~$ \|u_{\Lam}^{(1)}\|_{L^{4}(\R^{+};\dot H^{1})} $.
  
As mentioned in the introduction, any global solution belongs 
 to~$L^{4}(\R^{+};\dot H^{1})$. Thus, by definition of~$u_{\Lam}^{(1)}$, we have
 $$
 \|u_{\Lam}^{(1)}\| _{L^{4}(\R^{+};\dot H^{1})} \leq
 \sum_{J \in  \{1,\cdots,K\}}  \|u_{\Lam,J} \| _{L^{4}(\R^{+};\dot H^{1})}.
 $$
 Using a scaling argument, we infer 
 \begin{eqnarray}\label{estimateuLam1}
  \|u_{\Lam}^{(1)}\| _{L^{4}(\R^{+};\dot H^{1})}  & \leq &
  \sum_{J \in  \{1,\cdots,K\}}  \|u_{ J} \| _{L^{4}(\R^{+};\dot H^{1})} \nonumber\\
   & \leq  &  K \sup_{J \in \{1,\cdots,K\}}C_{J} \with \\
   C_{J}  & \eqdefa &   \|u_{J}   \|_{L^{\infty}(\R^{+};\dot H^{\frac12})} +   \|u_{J}   \|_{L^{2
 }(\R^{+};\dot H^{\frac32})} \nonumber.
 \end{eqnarray}

In order to estimate~$ \| F_{\Lam}\|_{L^{2}(\R^{+};\dot H^{-\frac12})} $, let
us  start by noticing that~$ F_{\Lam}$ is bounded uniformly in~$\Lam$  in the 
space~$L^{\frac43}(\R^{+};L^{2})$, by a constant
 depending on~$K$ and on the initial data. Indeed H\"{o}lder's inequality and Sobolev embeddings give
 \beno
 \| u_{\Lam,J}  \cdot \nabla u_{\Lam,J'}\|_{L^{\frac43}(\R^{+};L^{2})} &\leq&
   \| u_{\Lam,J}\|_{L^{4}(\R^{+}; L^{6})}
 \| \nabla u_{\Lam,J'}\|_{ L^{2}(\R^{+}; L^{3})}\\
  &\leq &     C   \| u_{\Lam, J}\|_{L^{4}(\R^{+}; \dot H^{1})}
 \| \nabla u_{ \Lam,J'}\|_{ L^{2}(\R^{+}; \dot H^{\frac12})},
 \eeno
 so that by scale invariance
 \beno
 \| F_{\Lam}\|_{L^{\frac43}(\R^{+};L^{2})}   &\leq &      C  \sum_{J \neq J'} \| u_{ J}\|_{L^{4}(\R^{+}; \dot H^{1})}
 \| \nabla u_{ J'}\|_{ L^{2}(\R^{+}; \dot H^{\frac12})} \\
 &\leq &      C K^{2} \sup_{J,J'}(C_{J} C_{J'}).
 \eeno
 So by interpolation  it is enough to prove that
  \beq\label{limFLam}
 \lim_{\Lam\rightarrow \infty}  \| F_{\Lam}\|_{L^{4}(\R^{+};\dot H^{-1})} = 0.
 \eeq
 Let~$J \neq J'$ be two integers in~$\{1,\dots , K\}$, and let~$\e > 0$ be given. 
 There exists a positive~$R$ and two vector fields~$\psi_{\e}$ and~$ \varphi_{\e}$
 in~${\mathcal D} (\R \times B(0,R))$ 
 such that
 $$
 \| \psi_{\e} -u_{J} \|_{L^{4}(\R^{+};\dot H^{1})}+
 \| \varphi_{\e} -u_{J'} \|_{L^{4}(\R^{+};\dot H^{1})} \leq \e.
 $$
 The support of~$T_{\Lam,J}\psi_{\e}$ (resp.~$T_{\Lam,J'}\vf_{\e}$) 
 is included in the  ball~$B(x_{J},R\Lam^{-1})$ (resp.~$B(x_{J'},R\Lam^{-1})$). 
 Thus we have
 \beq
 \label{limFLameq1}
 \Lam\geq 4 \delta^{-1}R\Longrightarrow T_{\Lam,J}\psi_{\e} \: 
 T_{\Lam,J'}\vf_{\e} =0.
 \eeq
 Then Sobolev embeddings as above give   the estimate
 $$\longformule
{\|u_{J} \otimes ( \varphi_{\e} -u_{J'})\|_{L^{4}(\R^{+};L^{2})} +
  \| ( \psi_{\e} -u_{J}) \otimes   \varphi_{\e} \|_{L^{4}(\R^{+};L^{2})}
}{  \leq C
\left(  \|u_{J}   \|_{L^{\infty}(\R^{+};\dot H^{\frac12})}
   \| \varphi_{\e} -u_{J'} \|_{L^{4}(\R^{+};\dot H^{1})}+
     \|\varphi_{\e}  \|_{L^{\infty}(\R^{+};\dot H^{\frac12})}
   \| \psi_{\e} -u_{J} \|_{L^{4}(\R^{+};\dot H^{1})}\right),}
 $$
 so that, using the scaling,
 \beq\label{CCjeps}
 \|u_{\Lam,J} \otimes ( T_{\Lam,J'}\vf_{\e} -u_{J'})\|_{L^{4}(\R^{+};L^{2})} +
  \| (T_{\Lam,J}\psi_{\e} -u_{J}) \otimes  T_{\Lam,J'}\vf_{\e} \|_{L^{4}(\R^{+};L^{2})} \leq C( C_{J} +
  C_{J'} ) \e.
  \eeq
Using~\refeq{limFLameq1},  it follows that for~$\Lam$
 large enough,
 $$
 \| F_{\Lam}\|_{L^{4}(\R^{+};\dot H^{-1})}  \leq C K^{2}  \e\sup_{J \in \{1,\cdots,K\}}C_{J},
$$
 and~(\ref{limFLam}) is proved.
 Plugging together that estimate with~(\ref{estimateuLam1}) gives~(\ref{conditionfacile}) for~ $\Lam$ large
 enough, and Proposition~\ref{stabildebilgeneral} is proved.
\end{proof}

\subsection{Proof of Theorems~\ref{stability} and~\ref{actionfractalonexample}}\label{plusdur}
Before starting the proofs, let us make a few comments on the transformation~$\tlk$ and
  state its
main properties. 
In all that follows, we shall consider only the action of~$\tlk$ on functions compactly supported in~$Q$.
First, one can  notice that if the family~$X$ of points satisfies~(\ref{distancedelta}), then
if~$\Lam \geq 4 \delta^{-1}$, 
$$
\mbox{\rm{supp}}\:\: T_\Lam^{J} f \:  \subset  \: Q _{\Lam }^{J} \eqdefa
\Bigl \{x \:  \Big/ \: d(x,x_{J}) \leq \Lambda^{-1}\Bigr\} \:  \subset  \:  Q _{ \delta}^{J}
\eqdefa  \Bigl\{x \:  \Big/ \: d(x,x_{J})
  \leq \frac14 \delta \Bigr\} .
$$
This implies immediately  that
\beq\label{lptlk}
\|\tlk f\|_{L^{p}} = \Lam^{ 1-\frac3p } K^{\frac1p}
\|f\|_{L^{p}} .
\eeq
 Then let us state   the following two lemmas, which are crucial for the proof of 
 Theorem~\ref{stability}
 and will be proved in Section~\ref{technicalities}.
\begin{lemma}\label{technicaltlamk}
{\sl
Let~$K \geq 1$ be an integer and~$\delta>0$ a real number.
There is  a constant~$C_{K,\delta} $  such that the following
results hold.  Let $r$ be in~$[1,\infty]$ and
consider a family~$X$ as in Definition~\ref{defTLamX}. Then for any real
number~$\Lambda$ in~$2^{\N}$
 greater than~$4 \delta^{-1}$ and for any~$ f \in {\mathcal D}(Q)$, we have
$$
 \| f\|_{\dot B^{-1}_{\infty,r}}-C_{K,\delta}\Lam^{-2 }  \|f\|_{\dot B^{-3}_{\infty,\infty}} \leq \|\tlk f\|_{\dot B^{-1}_{\infty,r}}
 \leq  \| f\|_{\dot B^{-1}_{\infty,r}} +C_{K,\delta} \Lam^{-2 }  \|f\|_{\dot B^{-3}_{\infty,\infty}}.
$$
Moreover the following estimate holds, where the constant~$C$ is universal:
\beq\label{h-sigma}
\|\tlk f\|_{\dot H^{-1}}
\leq C\sqrt K \Lam^{ - \frac32}\|f\|_{\dot H^{-1}} .
\eeq
}
\end{lemma}
{\bf Remark }  Let us point out that~$L^{1}$ is continuously included
in~$\dot B^{-3}_{\infty,\infty}$.

\begin{lemma}\label{heatheat}
{\sl Let~$K \geq 1$ be an integer and~$\delta>0$ a real number.
There is  a constant~$C_{K,\delta}$  such that the following
results hold.  Consider a family~$X$ as in Definition~\ref{defTLamX}. Then
 for any real number~$\Lambda$  in~$2^{\N}$
 greater than~$4 \delta^{-1}$ and for all divergence free vector fields~$f$ and~$g$
   in~${\mathcal D}(Q)$, we have
$$
\| {\bf P}(\St(t)\tlk f \cdot \nabla \St(t) \tlk g)\|_{E} \leq \| {\bf P}(\St(t) f \cdot \nabla \St(t) g)\|_{E}
+ C_{K,\delta} \Lam^{-3}\|f\|_{\dot H^{-1}}\|g\|_{\dot H^{-1} }.
$$
}
\end{lemma}

\subsubsection{End of the proof of Theorem~\ref{stability}}\label{endproof3}
Let us   consider a vector field~$u_{0} \in {\mathcal D}(Q)$ satisfying~(\ref{nonlineareps})
for
some~$\eta  \in ]0,1[$.       We know from Lemma~\ref{technicaltlamk}  that for any~ $ r \in [1,\infty]$
 and any~$\eta \in ]0,1[$,
    for any~ $\Lam$ greater than some~$\Lam_{0}$,  we have
\begin{equation}
\label{equivc-1}
\| u_{0} \|_{\dot B^{-1 }_{\infty,r}}
-\eta
\leq \|\tlk u_{0} \|_{\dot B^{-1 }_{\infty,r}} \leq \| u_{0} \|_{\dot B^{-1 }_{\infty,r}}
+ \eta.
\end{equation}

Next let us consider the smallness condition~(\ref{nonlinear}). By Lemma~\ref{heatheat} we know that
as soon as~$\Lam_{0}$ is large enough, then for any~$\Lam \geq \Lam_{0}$,
$$
\|{\bf P}( \St(t) \tlk u_{0}\cdot \nabla \St(t) \tlk u_{0})\|_{E} \leq \|{\bf P} (\St(t)  u_{0} \cdot \nabla \St(t)  u_{0})\|_{E}
+\eta.
$$
So we infer that
\beno
\| {\bf P}(\St(t) \tlk u_{0}\cdot \nabla \St(t) \tlk u_{0})\|_{E}
 &\leq&  C_{0}^{-1} \exp \left(
-C_{0}\bigl ( \|u_{0}\|_{\dot B^{-1}_{\infty,2}}+\eta\bigr) ^{4}
\right) \\
  &\leq&  C_{0}^{-1} \exp \left(
 -C_{0}  \| \tlk  u_{0}\|_{\dot B^{-1 }_{\infty,2}}^{4}  \right)
\eeno
due to~(\ref{equivc-1}). So Theorem~\ref{stability} is proved, up to the proof of Lemmas~\ref{technicaltlamk} and~\ref{heatheat}
which is the object of the coming section. \hfill $\square$

\subsubsection{The properties  of~$T_{\Lam,X}$}\label{technicalities}

In this    section, we are going to prove the properties  of
the transformation~$T_{\Lam,X}$ required in the proof of Theorem\refer{stability},
 namely Lemmas~\ref{technicaltlamk} and~\ref{heatheat}.
Before starting the proofs, let us give some more notation and prove preliminary results which will be used many times in the rest of this section.

We define
\beq
\label{definQtilde}
\widetilde Q_{\delta} = \bigcup_{J\in \{1,\dots,K\} }\widetilde Q_{\delta}^{J}, \quad \mbox{where}
\quad \widetilde Q_{\delta}^{J} \eqdefa \left\{
x \Big / d(x,Q_{\delta}^{J}) \leq \frac1{32} \delta
\right\},
\eeq
and we notice that this is a disjoint reunion.

\medbreak
The proof of Lemmas~\ref{technicaltlamk} and~\ref{heatheat} relies on the fact that
the Littlewood-Paley  theory is almost local. More precisely, let us recall Lemma 9.2.2
of\ccite{chemin14}.
\begin{lemma}
\label{lp2microlocfond}
{\sl
 For any positive integer~$N$ and any real number~$r$,  a
constant~$C_{N}$ exists such that the following result holds. Let~$F$ be a closed subset of~$\R^{3}$
and~$u$ a   distribution
in~$\dot B^{r}_{\infty,\infty}$   supported in~$F$; then for any couple~$(j,h)$ in~$\Z \times \R^+$
 such that~$2^{j}$ and~$2^{j} h$ are greater than~$1$, we have
$$
\norm {\Delta_{j} u}{L^\infty({}^cF_h)}\leq C_{N}2^{-jr}(2^{j} h)^{-N}
\|u\|_{\dot B^{r}_{\infty,\infty}},
$$
where~$F_h = \left\{x \in \R^{3} \: \Big / \: d(x,F ) \leq h\right\}$.
}
\end{lemma}

From this lemma, we deduce the following corollary.
\begin{corol}
\label{preliminarylemma}
{\sl Let~$K$, $\delta$ and~$X$ be as in Definition~\ref{defTLamX} and let~$M \in \N$ be
given. There  is a constant~$C_{ M}$  (depending only on~$M$)
  such that the following holds. For any~$\Lam \geq 4\delta^{-1}$,  for
any distribution~$f$ in~$\dot B^{-3}_{\infty,\infty}$, compactly supported  
in~$Q$ and
for any~$J \in \{ 1, \dots, K\}$, one has the following estimates:
\beq
\label{preliminarylemmaeq}
 \forall j \in \ZZ, \quad
\|\D_{j}T_\Lam^{J} f \|_{L^{\infty}( {}^{c} \widetilde Q_{\delta}^{J} )} \leq
C_{ M} \delta^{-(M+3)}\Lam^{-2}2^{-jM}  \|f\|_{\dot B^{-3}_{\infty,\infty}}.
\eeq
Moreover there is a universal constant~$C$ such that for any positive~$R$,
\beq
\label{claimsupports}
 \|f\|_{\dot B^{-1}_{\infty,r}} \leq
\Bigl\|\Bigl(2^{-j}\|\Delta_{j}f \|_{L^{\infty}( Q_{R})} \Bigr)_{j}
\Bigr\|_{\ell^r}
 + C R^{-2} \|f\|_{\dot B^{-3}_{\infty,\infty}},
 \eeq
where~$Q_{R} =  \left\{x \in \R^{3} \: \Big / \: d(x,Q) \leq R\right\}$.
}
\end{corol}

\begin{proof}
The first inequality is obvious when~$j$ is negative or when~$2^j \delta\leq 1$. Indeed we have the
scaling property
\beq \label{scaling}
\Delta_{j} \left(
f\bigl(\Lambda (\cdot - x_{J})\bigr)
\right) (x)= (\Delta_{j - \log_{2} \Lambda} f) (\Lambda (x  - x_{J})),
\eeq
so that for any~$s \in \R$,
$$
\|f (\Lambda (\cdot - x_{J}))\|_{\dot B^{s}_{\infty,\infty}} = \Lambda ^{s} \|f\|_{\dot B^{s}_{\infty,\infty}}
$$
Thus let us assume that~$2^j$ and~$2^j \delta$ are greater than~$1$. Using
Lemma\refer{lp2microlocfond}, we get
\beno
\|\D_{j}T_\Lam^{J} f \|_{L^{\infty}( {}^{c} \widetilde Q_{\delta}^{J} )}
 & \leq  &  C_{M}2^{3j} (2^{j} \delta)^{-M} \|T_\Lam^{J} f\|_{\dot B^{-3}_{\infty,\infty}}\\
 & \leq  & C_{M} 2^{-j(M-3)}  \delta^{-M} \Lam^{-2}\|f\|_{\dot B^{-3}_{\infty,\infty}}.
\eeno
In order to prove the second inequality, let us  note that, thanks to the triangle 
inequality and to the fact that~$\|\cdot\|_{\ell^r}\leq \|\cdot\|_{\ell^1}$, 
we have for any integer~$j_{0}$,
\[
\|f\|_{\dot B^{-1}_{\infty,r}} \leq
\Bigl\|\Bigl(2^{-j}\|\Delta_{j}f \|_{L^{\infty}( Q_{R})} \Bigr)_{j}
\Bigr\|_{\ell^r}
+ \sum_{j<j_{0}} 2^{-j}\|\D_{j}f\|_{L^\infty(\R^3)}
+\sum_{j\geq j_{0}}2^{-j}\|\D_{j} f\|_{L^\infty({}^c  Q_{R})}.
\]
Lemma\refer{lp2microlocfond} claims in particular that, if~$2^jR\geq 1$,
\[
\|\D_{j} f\|_{L^\infty({}^c Q_{R})}
\leq C R^{-3}\|f\|_{\dot B^{-3}_{\infty,\infty}}.
\]
Thus, if~$j_{0}$  is such that~$2^{j_{0}}R\geq 1$, we have, by definition of the norm
of~$\dot B^{-1}_{\infty,\infty}$,
\beno
\|f\|_{\dot B^{-1}_{\infty,r}} & \leq  &
\Bigl\|\Bigl(2^{-j}\|\Delta_{j}f \|_{L^{\infty}(  Q_{R})} \Bigr)_{j}
\Bigr\|_{\ell^r}
+\Bigl(\sum_{j<j_{0}}2^{2j}+
R^{-3}\sum_{j\geq j_{0}} 2^{-j}\Bigr) \|f\|_{\dot B^{-3}_{\infty,\infty}}\\
& \leq  &
\Bigl\|\Bigl(2^{-j}\|\Delta_{j}f \|_{L^{\infty}(  Q_{R})} \Bigr)_{j}
\Bigr\|_{\ell^r}
+ (2^{2j_{0}}+2^{-j_{0}}R^{-3})\|f\|_{\dot B^{-3}_{\infty,\infty}}.
\eeno
Choosing~$2^{j_{0}}\sim R^{-1}$ gives the result.
\end{proof}

\subsubsection{ Proof of Lemma~\ref{technicaltlamk}}
 We shall start by proving the second inequality, namely that
 \beq\label{secondinequality}
 \|\tlk f\|_{\dot B^{-1}_{\infty,r}} \leq \|f\|_{\dot B^{-1}_{\infty,r}} + C_{K,\delta} \Lam^{-2}
  \|f\|_{\dot B^{-3}_{\infty,\infty}} .
 \eeq
 Let us start  with   low
frequencies. We can  write that
\begin{eqnarray*}
\sum_{j < 0} 2^{-j} \|\Delta_{j} \tlk f\|_{L^\infty} & \leq &
\sum_{j < 0} 2^{-j}\sum_{J} \|\Delta_{j}T_\Lam^{J}  f\|_{L^\infty}  \\
 & \leq &  \sum_{J} \Bigl(\sum_{j< 0} 2^{2j}\Bigr)
 \|T_\Lam^{J} f\|_{ \dot B^{-3}_{\infty,\infty}}.
\end{eqnarray*}
Using the scaling equality~(\ref{scaling}) we get that
\beq
\label{lowfreq}
\sum_{j < 0} 2^{-j} \|\Delta_{j} \tlk f\|_{L^\infty}\leq
K\Lam^{-2}\|f\|_{\dot B^{-3}_{\infty,\infty}}.
\eeq
Now let us concentrate on the high frequencies. Recalling the definition given
in~(\ref{definQtilde}), let us start by considering the case when~$x \notin \widetilde Q_{\delta}$.
Using Inequality\refeq{preliminarylemmaeq} of Corollary\refer {preliminarylemma},
we can write (choosing~$M = 0$)
\beno
\|\D_{j}\tlk f\|_{L^\infty({}^{c}\wt Q_{\delta})} & \leq  &
\sum_{J} \|\D_{j}T_\Lam^{J} f\|_{L^\infty({}^{c}\wt Q_{\delta})}\\
& \leq & C K \delta^{-3} \Lam^{-2}\|f\|_{\dot B^{-3}_{\infty,\infty}}.
\eeno
Then we infer that
\beq
\label{prooflemmatlkeq111}
\sum_{j\geq 0}2^{-j}\|\D_{j}\tlk f\|_{L^\infty({}^{c}\wt Q_{\delta})}
\leq C K \delta^{-3} \Lam^{-2}\|f\|_{\dot B^{-3}_{\infty,\infty}}.
\eeq
Now let us consider the case when~$x \in \widetilde Q_{\delta}$. We can write
$$
\|\Delta_{j} \tlk f\|_{L^\infty(\widetilde Q_{\delta} )} \leq \sup_{J }
 \|\Delta_{j} \tlk f\|_{L^\infty(\widetilde Q^{J}_{\delta} )},
$$
and let us fix some~$J\in \{1,\dots,K \} $. We recall that
\[
\tlk f  = T_\Lam^{J}f  + \sum_{J' \neq J} T_\Lam^{J'}f ,
\]
and let us start with the estimate of~$\tlk f - T_\Lam^{J}f $. We have
\[
\|\D_{j}(\tlk f-T_\Lam^{J} f)\|_{L^\infty(\wt Q_{\delta}^J)} \leq
\sum_{J'\not = J} \|\D_{j}T_\Lam^{J'} f\|_{L^\infty({}^c\wt Q_{\delta}^{J'})} .
\]
Using Inequality\refeq{preliminarylemmaeq} of Corollary\refer {preliminarylemma},
we get that
\[
\|\D_{j}(\tlk f-T_\Lam^{J} f)\|_{L^\infty({}^c\wt Q_{\delta}^J)} \leq
C K \delta^{-3} \Lam^{-2}\|f\|_{\dot B^{-3}_{\infty,\infty}}.
\]
Thus we infer
\beq
\label{highfreq1}
\sum_{j\geq 0} 2^{-j}\|\D_{j}(\tlk f-T_\Lam^{J} f)\|_{L^\infty(\wt Q_{\delta}^J)}
\leq C K \delta^{-3} \Lam^{-2}\|f\|_{\dot B^{-3}_{\infty,\infty}}.
\eeq
Now let us examine the term~$\|T_\Lam^{J}f\|_{L^\infty(\wt Q_{\delta}^J)}$.  From~(\ref{scaling}) we get
\beno
 \Bigl\|\Bigl({\bf 1}_{j \geq 0}
2^{-j} \|\Delta_{j} T_\Lam^{J}f\|_{L^\infty(\widetilde Q^{J}_{ \delta} )}\Bigr)_{j}
\Bigr\|_{\ell^r} & \leq & \Lam
\Bigl\|\Bigl({\bf 1}_{j \geq 0}
2^{-j} \|\Delta_{j-  \log_{2} \Lam}f\|_{L^\infty(\R^{3})}\Bigr)_{j}
\Bigr\|_{\ell^r}\nonumber \\
& \leq &  \|f\|_{\dot B^{-1}_{\infty,r}}.
\eeno
Once noticed that~$\|\cdot\|_{\ell^r}\leq \|\cdot\|_{\ell^1}$, we plug
together that estimate with~(\ref{lowfreq}), (\ref{prooflemmatlkeq111}) and~(\ref{highfreq1})
to conclude the proof of~(\ref{secondinequality}).

Let us bound from below~$\|\tlk f\|_{\dot B^{-1}_{\infty,r}}$.
As~$\|g\|_{L^{\infty}(\wt Q_{\delta})}=
\ds\sup_{J} \|g\|_{L^{\infty}(\wt Q^{J}_{\delta})}$, we have
\beno
 \|\tlk f\|_{\dot B^{-1}_{\infty,r}} &\geq &
\Bigl\|\Bigl(2^{-j}
\|\Delta_{j} \tlk f \|_{L^{\infty}(\wt Q_{\delta})}
\Bigr)_{j}\Bigr\|_{\ell^r}\\
 & \geq & \Bigl\|\Bigl(2^{-j}\sup_{J}
\|\Delta_{j} \tlk f \|_{L^{\infty}(\wt Q^{J}_{ \delta})}
\Bigr)_{j}\Bigr\|_{\ell^r}\\
  &\geq & \Bigl\|\Bigl(2^{-j}
\|\Delta_{j} \tlk f \|_{L^{\infty}(\wt Q^{J_{0}}_{\delta})}
\Bigr)_{j}\Bigr\|_{\ell^r} \eeno
for some~$J_{0}$ in~$\{1,\dots,K \} $. Using the fact
that~$\|\cdot\|_{\ell^r}\leq \|\cdot\|_{\ell^1}$, we can write that
\[
\longformule{
 \|\tlk f\|_{\dot B^{-1}_{\infty,r}} \geq
 \Bigl\|\Bigl(2^{-j}
\|\Delta_{j} T_\Lam^{J_{0}} f \|_{L^{\infty}(\wt Q^{J_{0}}_{\delta})}
\Bigr)_{j}\Bigr\|_{\ell^r}
}
{
{}- \sum_{j< 0}2^{-j} \Bigl\|
\Delta_{j}T_\Lam ^{J_{0}} f
\Bigr \|_{L^{\infty}(\wt Q^{J_{0}}_{\delta})}
- \sum_{j\geq 0}2^{-j} \Bigl\|
\Delta_{j} (\tlk f-T_\Lam^{J_{0}} f)
\Bigr \|_{L^{\infty}(\wt Q^{J_{0}}_{\delta})}.
}
\]
Using\refeq{lowfreq} and\refeq{highfreq1}, we infer that
\beq
\label{eqproofminorC-122}
 \|\tlk f\|_{\dot B^{-1}_{\infty,r}} \geq
 \Bigl\|\Bigl(2^{-j}
\|\Delta_{j}T_\Lam^{J_{0}}  f \|_{L^{\infty}(\wt Q^{J_{0}}_{\delta})}
\Bigr)_{j}\Bigr\|_{\ell^r}-C_{K,\delta}\Lam^{-2}\|f\|_{\dot B^{-3}_{\infty,\infty}}.
\eeq
By scaling and translation, we have
\[
 \Bigl\|\Bigl(2^{-j}
\|\Delta_{j}T_\Lam^{J_{0}}  f \|_{L^{\infty}(\wt Q^{J_{0}}_{k})}
\Bigr)_{j}\Bigr\|_{\ell^r} =
 \Bigl\|\Bigl(2^{-j}
\|\Delta_{j} f \|_{L^{\infty}(Q_{\Lam,\delta})}
\Bigr)_{j}\Bigr\|_{\ell^r}
\]
where~$Q_{\Lam,\delta}$ is   the cube  of size~$2 \Lam  \delta$.
Using\refeq{claimsupports} with~$R=2 \delta \Lam$ and\refeq{eqproofminorC-122}, we infer that
\[
\|\tlk f\|_{\dot B^{-1}_{\infty,r}} \geq  \|f\|_{\dot B^{-1}_{\infty,r}}  -
C_{K,\delta}\Lam^{-2}\|f\|_{\dot B^{-3}_{\infty,\infty}}.
\]
This concludes the proof of the first part of the lemma.

\medskip

\noindent Now let us prove the second part of the lemma, namely Estimate~(\ref{h-sigma}) on the~$\dot H^{-1}$ norm.  Let~$f \in \cD(Q)$ be given. Stating~$f_{m}\eqdefa
-\cF^{-1}(i\xi_{m}|\xi|^{-2}\wh f)$, we can write
$$
f = \sum_{m = 1}^{3} \partial_{m} f_{m}, \with \|f\|_{\dot H^{-1 }}
 \sim \sum_{m = 1}^{3}  \|f_{m}\|_{L^{2}}.
$$
Let us recall that
$$
\|\tlk f\|_{\dot H^{-1 }} = \supetage{g \in {\cD} (Q)}{\|g\|_{\dot H^{ 1 }} \leq 1} \int_{\R^{3}}
\tlk f (x) g(x) \: dx.
$$
Let~$\chi \in {\cD} (Q)$ be equal to one on the support of~$g$. We have
\beno
 \int_{\R^{3}}\tlk f (x) g(x) \: dx &=&\Lam \sum_{J} \int_{\R^{3}}  f (\Lam (x-x_{J})) g(x) \: dx \\
   &=&\Lam^{-2 }\sum_{J}\sum_{m}\int_{\R^{3}}    \partial_{m} f_{m}  (x) g(\Lam^{-1}x+x_{J})\: dx ,
\eeno
so after an integration by parts and a change of variables again, we infer that
\beno
 \int_{\R^{3}}\tlk f (x) g(x) \: dx  & = &
 -\Lam^{-3 }\sum_{J}\sum_{m}\int_{\R^{3}}   \chi (x) f_{m}  (x) (\partial_{m}g)(\Lam^{-1}x+x_{J})\: dx \\
  & = & -\Lam^{-1} \sum_{m}\int_{\R^{3}} \tlk(  \chi  f_{m})  (x)   \partial_{m}g (x) \: dx.
\eeno
In particular we get that
\beno
\|\tlk f \|_{\dot H^{-1 }}  & \leq & C \Lam^{-1}
\sum_{m} \|  \tlk(  \chi  f_{m})\|_{L^{2}} \\
  &\leq& C\Lam^{-1} \sum_{m} \|   f_{m} \|_{L^{2}}  \Lam^{-\frac 12}
  \sqrt K \\
    &\leq& C  \Lam^{-\frac {3}2}   \sqrt K\|f\|_{\dot H^{-1 }}  ,
\eeno
and the result is proved.
\hfill $\square$

\subsubsection{ Proof of Lemma\refer{heatheat}}   First, we observe that
\beq
\label{div0lf}
\St(t)\tlk f \cdot \nabla \St(t) \tlk g =
\sum_{\ell=1}^3 \partial_{\ell}\Bigl(\St(t)\tlk f^\ell \St(t)  \tlk g\Bigr),
\eeq
so  using Bernstein's inequalities, we can write
\beno
\cE_{j} & \eqdefa & \Bigl\|\Delta_{j}\Bigl(\St(t)\tlk f
 \cdot \nabla \St(t) \tlk g \Bigr)\Bigr\|_{L^{1}(\R^{+};L^{\infty})}\\
 &&\qquad\qquad\qquad\qquad{}+
  \Bigl\|t^\frac 1 2\Delta_{j}\Bigl(\St(t)\tlk f
 \cdot \nabla \St(t) \tlk g \Bigr)\Bigr\|_{L^{2}(\R^{+};L^{\infty})}\\
  &\leq & C2^{4j}\biggl(\|\St(t)\tlk f\|_{L^{2}(\R^{+};L^{2})}
\|\St(t) \tlk g\|_{L^{2}(\R^{+};L^{2})}\biggr.\\
&& \biggl.\qquad\qquad\qquad\qquad{} +
\|t^\frac 1 2\St(t)\tlk f\|_{L^{\infty}(\R^{+};L^{2})}
\|\St(t) \tlk  g\|_{L^{2}(\R^{+};L^2)}\biggr)\\
& \leq & C 2^{4j}
\| \tlk f\|_{\dot H^{-1}}\| \tlk g\|_{\dot H^{-1}}.
\eeno
Using Lemma\refer{technicaltlamk}, we get
$$
\cE_{j}\leq C 2^{4j} K \Lam^{-3 }
\|f\|_{\dot H^{-1}}\|g\|_{ \dot H^{-1}  }.
$$
We therefore infer a bound on the low frequencies:
\beq
\label{proofheattheateq1}
\sum_{j \leq 0} 2^{-j} \cE_{j}\leq C   K \Lam^{-3 }
\|f\|_{\dot H^{-1}}\|g\|_{ \dot H^{-1} }.
\eeq
The high frequencies are   more delicate to estimate. Let us write that
\beno
\St(t)\tlk f \cdot \nabla \St(t)\tlk g & = & H_{\Lam,X}(t)+K_{\Lam,X}(t)\with\\
H_{\Lam,X} (t)& \eqdefa & \sum_{J\not =J'} \sum_{\ell=1}^3
\partial_{\ell}\Bigr(\St (t)T_\Lam^{J} f^\ell \St(t) T_\Lam^{J'}g\Bigr)
\andf\\
K_{\Lam,X} (t)& \eqdefa & \sum_{J} \sum_{\ell=1}^3
\partial_{\ell}\Bigl
(\St (t) T_\Lam^{J}f^\ell\St(t) T_\Lam^{J}g\Bigr).
\eeno
We observe that
\beno
B^{J,J'}_{\Lam, X}(f,g) &\eqdefa &
\partial_{\ell}\Bigr(\St (t)T_\Lam^{J} f^\ell \St(t)T_\Lam^{J'}g\Bigr)\\
 & = &
 \frac 1 {(4\pi t)^{3}}\int_{\R^6}
\partial_{x_{\ell}}\exp \Bigl(-\frac {|x-y|^2+|x-z|^2}  {4t} \Bigr)
T_\Lam^{J}f^\ell (y)T_\Lam^{J'}g(z)dydz\\
& = & \frac {-1} {(4\pi t)^{3}}\int_{\R^6}
\frac {2x_{\ell}-y_{\ell}-z_{\ell}}{2t}
\exp \Bigl(-\frac {|x-y|^2+|x-z|^2}  {4t} \Bigr) T_\Lam^{J}f^\ell (y)T_\Lam^{J'}g(z)dydz .
\eeno
Due to the distance between~$x_{J} $ and~$x_{J'}$,
one gets that a smooth bounded function (as well as all its derivatives)~$\chi$
on~$\R$ exists such that~$\chi$ vanishes identically near $0$ and such that
\beno
B^{J,J'}_{\Lam, X}(f,g) (t,x) &= & \frac1{(4\pi t)^{3}}\int_{\R^6}\Theta_{\delta}(t,x,y,z)
T_\Lam^{J}f (y)T_\Lam^{J'}g(z)dydz \with\\
\Theta_{\delta}(t,x,y,z) & \eqdefa & \frac 1 {t^{\frac  1 2}} \chi\left(\frac{|x-y|^2+|x-z|^2}{C\delta^{2}}\right)
\frac {2x_{\ell}-y_{\ell}-z_{\ell}}{2t^{\frac 1 2}}
\exp \Bigl(-\frac {|x-y|^2+|x-z|^2}  {4t} \Bigr) \cdotp
\eeno
As we have
\beq
\label{inegtensorproductHneg}
\|a\otimes b\|_{\dot H^{-2}(\R^6)}\leq
\|a\|_{\dot H^{-1}(\R^3)}\|b\|_{\dot H^{-1}(\R^3)},
\eeq
 we infer, using the scaling, that
\beno
\|B^{J,J'}_{\Lam, X}(f,g) (t,\cdot)\|_{L^\infty} & \leq &\frac1{(4\pi t)^{3}}
\sup_{x\in \R^3}\|\Theta_{\delta}(t,x,\cdot)\|_{\dot H^2(\R^6)}
 \|T_\Lam^{J}f \|_{\dot H^{-1}(\R^3)}\|T_\Lam^{J'}g\|_{\dot H^{-1}(\R^3)}\\
  & \leq  & \frac1{(4\pi t)^{3}}
  \sup_{x\in \R^3}\|\Theta_{\delta}(t,x,\cdot)\|_{\dot H^2(\R^6)}\Lam^{-3}
 \|f \|_{\dot H^{-1}(\R^3)}\|g\|_{\dot H^{-1}(\R^3)}.
\eeno
It is obvious that
\[
|\nabla ^2\Theta_{\delta}(t,x,y,z)|\leq  
\frac  C {t^{\frac  3 2}} e^{-\frac \delta {C t}}
\exp \Bigl(-\frac {|x-y|^2+|x-z|^2}  {8t} \Bigr)
\]
and thus that
\[
\|B^{J,J'}_{\Lam, X}(f,g) (t,\cdot)\|_{L^\infty} \leq
\frac C{t^3} \Lam^{-3}e^{-\frac \delta {Ct}} \|f \|_{\dot H^{-1}(\R^3)}\|g\|_{\dot H^{-1}(\R^3)}.
\]
We immediately infer, since~$\|\Delta_{j} \P a\|_{L^{\infty}} \leq C \|\Delta_{j} a\|_{L^{\infty}}$,  that
\beq
\label{Div0formule}
 \sum_{j\geq 0}2^{-j}\Bigl(\|\D_{j}{\bf P}H_{\Lam,X}\|_{L^1(L^\infty)}
 +\|t^{\frac 1 2}\D_{j}{\bf P}H_{\Lam,X}\|_{L^2(L^\infty)}\Bigr)
 \leq  C_{\delta}\Lam^{-3} \|f \|_{\dot H^{-1}}\|g\|_{\dot H^{-1}}.
 \eeq
Now let us consider the term~$K_{\Lam,X}$. To start with, let us write  
\beno
\|\D_{j}{\bf P}K_{\Lam,X}(t,\cdot)\|_{L^\infty(\R^3)} & \leq &
\|\D_{j}{\bf P}K_{\Lam,X}(t,\cdot)\|_{L^\infty(\wt Q_{\delta})}
+\|\D_{j}{\bf P}K_{\Lam,X}(t,\cdot)\|_{L^\infty({}^c\wt Q_{\delta})}\\
& \leq  & \sup_{J}\|\D_{j}{\bf P}K_{\Lam,X}(t,\cdot)\|_{L^\infty(\wt Q_{\delta}^J)}
+\|\D_{j}{\bf P}K_{\Lam,X}(t,\cdot)\|_{L^\infty({}^c\wt Q_{\delta})}.
\eeno
By definition of~$K_{\Lam,X}$, and denoting~$\wt \Delta_{j} = \Delta_{j} \P$,  we get
\beno
\|\D_{j}{\bf P}K_{\Lam,X}(t,\cdot)\|_{L^\infty(\R^3)} & \leq &
\sup_{J}\|\D_{j}{\bf P}(\St(t)T_\Lam^{J}f\nabla \St(t)T_\Lam^{J}g)\|_{L^\infty(\R^3)}\\
&&
\qquad\quad{}+ \sup_{J}\Bigl\|\wt \D_{j}
\sum_{J'\not =J}
\partial_{\ell}\Bigr(\St (t)T_\Lam^{J'} f^\ell \St(t) T_\Lam^{J'}g\Bigr)
\Bigr\|_{L^\infty(\wt Q_{\delta}^J)}\\
 &  &
\qquad\quad\qquad{}+ \sum_{J'}\Bigl\|\wt \D_{j}
\partial_{\ell}\Bigr(\St (t)T_\Lam^{J'} f^\ell \St(t)T_\Lam^{J'}g\Bigr)
\Bigr\|_{L^\infty({}^c\wt Q_{\delta})}\\
&\leq &
\sup_{J}\|\D_{j}{\bf P}(\St(t)T_\Lam^{J}f\nabla \St(t)T_\Lam^{J}g)\|_{L^\infty(\R^3)}\\
&&
\qquad\quad{}+ \sup_{J}
\sum_{J'\not =J}\Bigl\|\wt \D_{j}
\partial_{\ell}\Bigr(\St (t)T_\Lam^{J'} f^\ell \St(t) T_\Lam^{J'}g\Bigr)
\Bigr\|_{L^\infty(\wt Q_{\delta}^J)}\\
 &  &
\qquad\quad\qquad{}+ \sum_{J'}\Bigl\|\wt \D_{j}
\partial_{\ell}\Bigr(\St (t)T_\Lam^{J'} f^\ell \St(t)T_\Lam^{J'}g\Bigr)
\Bigr\|_{L^\infty({}^c\wt Q_{\delta})}\\
&\leq &
\sup_{J}\|\D_{j}{\bf P}(\St(t)T_\Lam^{J}f\nabla \St(t)T_\Lam^{J}g)\|_{L^\infty(\R^3)}\\
 &  &
\qquad\quad\qquad{}+ \sum_{J'}\Bigl\|\wt \D_{j}
\partial_{\ell}\Bigr(\St (t)T_\Lam^{J'} f^\ell \St(t) T_\Lam^{J'}g\Bigr)
\Bigr\|_{L^\infty({}^c\wt Q^{J'}_{\delta})}.
\eeno
By translation and scaling we infer that
\[
\longformule{
\|\D_{j}{\bf P}K_{\Lam,X}(t,\cdot)\|_{L^\infty(\R^3)}  \leq \Lam\|\D_{j-\log_{2}
\Lam}{\bf P}(\St(t)f\nabla \St(t)g)\|_{L^\infty(\R^3)}
}
{
{}+   \sum_{J'}\Bigl\|\wt \D_{j}
\partial_{\ell}\Bigr(\St (t)T_\Lam^{J'} f^\ell \St(t)T_\Lam^{J'}g\Bigr)
\Bigr\|_{L^\infty({}^c\wt Q^J_{\delta})}.
}
\]
By definition of~$\St(t)$, we have, for some~$\wt h \in {\mathcal S}(\R^{3})$, 
\beno
B^{J'}_{\Lam,j}(f,g)(t,x) & \eqdefa & \wt \Delta_{j}
\partial_{\ell}\Bigr(\St (t)T_\Lam^{J'} f^\ell \St(t)T_\Lam^{J'}g\Bigr)(t,x)\\
& = &\frac {2^{3j}}{(4\pi t)^3}\int_{\R^9}\wt h(2^j(x-x'))
\partial_{x'_{\ell}}\exp \Bigl(-\frac {|x'-y|^2+|x'-z|^2}  {4t} \Bigr)\\
&&\qquad\qquad\qquad\qquad\qquad\qquad\qquad\qquad{}\times
T_\Lam^{J'}f(y)T_\Lam^{J'}g(z)dx'dydz.
\eeno
Now if~$x$ is in~${}^c\wt Q_{\delta}^{J'}$ and~$y$ in~$\wt Q_{\delta}^{J'}$,
one gets that a smooth bounded function (as well as all its derivatives)~$\chi$
on~$\R$ exists such that~$\chi$ vanishes identically near $0$ and has value~$1$
outside a ball centered at the origin, and such that
\beno
B^{J'}_{\Lam,j}(f,g) & = & B^{J',1}_{\Lam,j}(f,g)+B^{J',2}_{\Lam,j}(f,g)\with\\
 B^{J',1}_{\Lam,j}(f,g) (t,x) & \eqdefa &
 \frac {2^{3j}}{(4\pi t)^3}\int_{\R^9}\wt h(2^j(x-x'))\chi( \frac{ |x'-x|^2}{C\delta^{2}})
\\
&&\qquad{}\times
\partial_{x'_{\ell}}\exp \Bigl(-\frac {|x'-y|^2+|x'-z|^2}  {4t} \Bigr)
T_\Lam^{J'}f(y)T_\Lam^{J'}g(z)dx'dydz\andf\\
 B^{J',2}_{\Lam,j}(f,g)(t,x)&\eqdefa &
 \frac {2^{3j}}{(4\pi t)^3}\int_{\R^9}\wt h(2^j(x-x'))\chi( \frac{ |x'-y|^2}{C\delta^{2}})
\\
&&\qquad{}\times
\partial_{x'_{\ell}}\exp \Bigl(-\frac {|x'-y|^2+|x'-z|^2}  {4t} \Bigr)
T_\Lam^{J'}f(y)T_\Lam^{J'}g(z)dx'dydz.
\eeno
By integration by parts,  we get that
\[
 B^{J',1}_{\Lam,j}(f,g) (t,x) = - 2^{3j}\!\!\int_{\R^9}\!\!\partial_{x'_{\ell}}
\Bigl (\wt h(2^j(x-x'))\chi( \frac{ |x'-x|^2}{C\delta^{2}})\Bigr)
(\St(t)T_\Lam^{J'}f)(x')(\St(t)T_\Lam^{J'}g)(x')dx'.
\]
By the Leibnitz formula  we have,
\[
\longformule {2^{3j}\partial_{x'_{\ell}}
\Bigl (\wt h(2^j(x-x'))\chi( \frac{ |x'-x|^2}{C\delta^{2}})\Bigr)  = 
2^{4j}(\partial_{x_{\ell}}\wt h)(2^{j}(x-x'))\chi( \frac{ |x'-x|^2}{C\delta^{2}})}
{
{}+C_{\delta}2^{3j}\wt h(2^{j}(x-x'))  (x_{\ell}-x'_{\ell})\chi( \frac{ |x'-x|^2}{C\delta^{2}})\\
}
\]
Using the properties of the function~$\chi$, we infer 
\[
\Bigl|2^{3j}\partial_{x'_{\ell}}
\Bigl (\wt h(2^j(x-x'))\chi( \frac{ |x'-x|^2}{C\delta^{2}})\Bigr) \Bigr|
\leq  C_\delta \underline h(2^{j}(x-x'))
\]
for some bounded  function~$\underline h$. Thus, by integration, we infer that
\[
\| B^{J',1}_{\Lam,j}(f,g)(t,\cdot)\|_{L^\infty}
\leq C_{\delta} \|\St(t)T_\Lam^{J'}f\|_{L^2}\|\St(t)T_\Lam^{J'}g\|_{L^{2}}.
\]
By definition of Besov spaces, we deduce that
\[
\| B^{J',1}_{\Lam,j}(f,g)\|_{L^1(L^\infty)}
+\|t^{\frac 1 2} B^{J',1}_{\Lam,j}(f,g)\|_{L^2(L^\infty)}
\leq C_{\delta} \|T_\Lam^{J'}f\|_{\dot H^{-1}}\|T_\Lam^{J'}g\|_{\dot H^{-1}}.
\]
By scaling, we infer that
\beq
\label{eqdemoestimbilieanhf1}
\| B^{J',1}_{\Lam,j}(f,g)\|_{L^1(L^\infty)}
+\|t^{\frac 1 2} B^{J',1}_{\Lam,j}(f,g)\|_{L^2(L^\infty)}
\leq C_{\delta} \Lam^{-3}\|f\|_{\dot H^{-1}}\|g\|_{\dot H^{-1}}.
\eeq
In order to estimate~$B^{J',2}_{\Lam,j}(f,g)$, let us write
\beno
B^{J',2}_{\Lam,j}(f,g)(t,x)  & =   &\frac {1}{(4\pi t)^3} \int_{\R^6}\Theta_{\delta,j}(t,x,y,z)
T_\Lam^{J'}f(y)T_\Lam^{J'}(z)dydz\with\\
\Theta_{\delta,j}(t,x,y,z)  & \eqdefa  &\frac {2^{3j}}{t^{\frac 1 2}}
\int_{\R^3}\wt h(2^j(x-x'))\chi( \frac{ |x'-y|^2}{C\delta^{2}})\\
&&\qquad\qquad{}\times
\frac {2x'_{\ell}-y_{\ell}-z_{\ell}}{2t^{\frac 1 2}}
\exp \Bigl(-\frac {|x'-y|^2+|x'-z|^2}  {4t}\Bigr) dx'.
\eeno
Using\refeq{inegtensorproductHneg},  the definition of the Besov norm and the
scaling property, we deduce that
\beno
\|B^{J',2}_{\Lam,j}(f,g)(t,\cdot)\|_{L^\infty} & \leq  & \sup_{x\in \R^3}
\|\nabla_{y,z}^2 \Theta_{\delta,j}(t,x,\cdot,\cdot)\|_{L^2(\R^6)}
\|T_\Lam^{J'}f\|_{\dot H^{-1}}\|T_\Lam^{J'}g\|_{\dot H^{-1}}\\
& \leq & \sup_{x\in \R^3}
\|\nabla_{y,z}^2 \Theta_{\delta,j}(t,x,\cdot,\cdot)\|_{L^2(\R^6)}
\Lam^{-3}\|f\|_{\dot H^{-1}}\|g\|_{\dot H^{-1}}.
\eeno
A straightforward computation shows that
\[
\sup_{x\in \R^3}
\|\nabla_{y,z}^2 \Theta_{\delta,j}(t,x\cdot,\cdot)\|_{L^2(\R^6)}
\leq  \frac { C_{\delta}}  {t^{\frac 3 2}} e^{- \frac \delta {Ct}}.
\]
Thus, we get that
\[
\| B^{J',2}_{\Lam,j}(f,g)\|_{L^1(L^\infty)}
+\|t^{\frac 1 2} B^{J',2}_{\Lam,j}(f,g)\|_{L^2(L^\infty)}
\leq
C_{\delta}\Lam^{-3}\|f\|_{\dot H^{-1}}\|g\|_{\dot H^{-1}}
\]
Using\refeq{proofheattheateq1},\refeq{Div0formule} and\refeq{eqdemoestimbilieanhf1} we infer that
\beno
\|{\bf P}(\St(t)\tlk f\nabla \St(t)\tlk g)\|_{E} & \leq &
\sum_{j}2^{-j}\Lam \Bigl
(\|\D_{j-\log_{2\Lam}}{\bf P}\St(t)f\nabla\St(t)g\|_{L^1(L^\infty)}\\
&& \qquad\qquad\quad
{}+\|t^{\frac 1 2}\D_{j-\log_{2\Lam}}{\bf P}\St(t)f\nabla\St(t)g
\|_{L^2(L^\infty)}\Bigr)\\
&& \qquad\qquad\qquad\qquad\qquad\qquad\quad
{}+C_{K,\delta}\Lam^{-3}\|f\|_{   \dot H^{-1}}\|g\|_{\dot H^{-1}}\\
& \leq & \|{\bf P}(\St(t) f\nabla \St(t) g)\|_{E}
+ C_{K,\delta}\Lam^{-3}\|f\|_{\dot H^{-1}}\|g\|_{\dot H^{-1}}.
\eeno
That ends the proof of Lemma~\ref{heatheat}. \hfill $\square$

\subsubsection{Proof of Theorem~\ref{actionfractalonexample}} The proof is 
straightforward: in order to apply
Theorem~\ref{stability}, we define~$\eta>0$ and we need to find~$\Lam_{0}$ uniform in~$\e$ so that,
according to Lemmas~\ref{technicaltlamk} and~\ref{heatheat}, the following two conditions are satisfied:
$$
\Lam_{0}^{-3} C_{\delta,K} \|u_{0,\varepsilon}\|_{\dot H^{-1}}^{2} = \eta \quad \mbox{and} \quad
\Lam_{0}^{-2} C_{\delta,K} \|u_{0,\varepsilon}\|_{\dot B^{-3}_{\infty,\infty}}  = \eta.
$$
Due to Corollary~\ref{cordonneeinitiale} this is trivially possible as soon as~$\alpha>0$. \hfill$ \Box$


\end{document}